\documentstyle{jtbart}

\makeatletter
\def\@begintheorem#1#2{\it \trivlist \item[\hskip
\labelsep{\bf #2\ #1.}]}
\def\@opargbegintheorem#1#2#3{\it \trivlist
      \item[\hskip \labelsep{\bf #1\ #2\ (#3).}]}
\makeatother

\newtheorem{them}{Theorem}[section]
\newtheorem{quest}{Question}
\newtheorem{lemm}[them]{Lemma}
\newtheorem{cor}[them]{Corollary}
\newtheorem{prop}[them]{Proposition}
\newtheorem{thm}[them]{Theorem}

\newcommand{\jtbnumpar}[1]{\refstepcounter{them}
\trivlist
\item[\hskip \labelsep{\bf \thethm \ #1.}]}

\newcommand{\jtbdef}{\jtbnumpar{Definition}}
\def\jtbnot{\jtbnumpar{Notation}}
\let\endjtbdef=\par

\newtheorem{ssnote}{COMMENT}

\def\PROOF #1.{\par\noindent{\it Proof#1}.\ \ignorespaces}
\def\proof #1.{\par\noindent{\it Proof#1}.\ \ignorespaces}

\def\subm{\leq}

\def\Mscr{{\cal M}}

\let\?=\joinrel

\let\leftv=^
\let\rightv=^

\def\endproof{\vskip 5pt plus1pt minus1pt }

\def\id{\mathop{\rm ID}}
\def\fil{\mathop{\rm FIL}}

\def\cub{\mathop{\rm cub}}
\def\inv{\mathop{\rm inv}}
\def\som{\mathop{{\rm sm}_{1}}}
\def\sm{\mathop{\rm sm}}
\def\acc{\mathop{\rm acc}}
\def\nacc{\mathop{\rm nacc}}

\def\k/{\kern.2em}    

\def\sm{\mathop{\rm sm}}

\def\otp{\mathop{\rm otp}}

\def\sing{\mathop{\rm sing}}
\def\cpr{\mathop{\rm cpr}}

\def\cf{\mathop{\rm cf}}

\def\sup{\mathop{\rm sup}}

\def\lg{\mathop{\rm lg}}

\def\REI{{\rm REI}}

\def\ls{\rm LS}               %

\let\bet=\Im        %

 \mathchardef\bet="0B69

\def\down{\smash{\mathchar"0223}}

\let\union=\cup             %

\edef\bigcup{\mathop{\textstyle\mathchar\the\bigcup}}

\let\inter=\cap             %

\edef\bigcap{\mathop{\textstyle\mathchar\the\bigcap}}

\edef\bigwedge{\mathop{\textstyle\mathchar\the\bigwedge}}

\edef\bigvee{\mathop{\textstyle\mathchar\the\bigvee}}

\edef\sum{\mathop{\textstyle\mathchar\the\sum}}
\def\ind #1#2#3{#1 \mathbin{\mathop{\down}_{#2}} #3}
\def\math&{\ \& \ }

\def\force {\mathrel^\joinrel\rightarrow}
\def\force {\mathrel{\scriptstyle\mathrel^\joinrel\rightarrow}}
\def\forceq {\mathrel{\mathop{\force}\limits_{\textstyle\texsim}}}
\def\forceq{\mathrel^\joinrel
 \mathrel{\mathop{\rightarrow}\limits_{\smash{\textstyle\texsim}}}}
\def\forceq{\mathrel{\scriptstyle\mathrel^\joinrel
 \mathrel{\smash{\mathop{\rightarrow}\limits_{\smash{\raise
 2pt\hbox{$\scriptstyle\texsim$}}}}}}}

\let\exclaim=!                 %
\let\iso\approx             %
\let\texsim=\sim         %
\let\conj\sim             %
\def\conjp #1 {\conj_{#1}}     %
\let\sim\simeq             %
\def\0bar{\bar 0}         %
\def\1bar{\bar 1}         %
\def\abar{\overline a}

\def\cbar{\overline c}

\def\hbar{\overline h}

\def\Mun{\underline M}
\def\Nun{\underline N}
\def\Lun{\underline L}

\let\ocirc\circ    

\def\tp{\rm tp}

\def\Proof{Proof}

\def\bet{]}

\def\Pun{\underline P}

\def\bet{]}
\def\restrict{ | }

\def\Mun{{\underline M}}
\def\Nun{{\underline N}}
\def\Dscr{{\cal D}}
\def\kd{{\cal D}}

\def\chik{\chi_{{\bf K}}}

\def\text#1{\ifmmode\leavevmode\hbox{#1}\else
   \typeout{Warning: \string\text \space used outside math mode!}
   \begingroup\hbox{#1}\endgroup\fi}

\def\chin{\chi_{1}({\bf K})}
\title{The Primal Framework II: Smoothness}
\author{J. Baldwin
\thanks{Partially supported by N.S.F. grant 8602558 }\\
Department of Mathematics\\
University of Illinois, Chicago
\and
S. Shelah
Department of Mathematics\\
Hebrew University of Jerusalem
\thanks{Both authors thank Rutgers University and the U.S. Israel
Binational Science foundation for their support of this project.
This is item 360 in Shelah's bibliography.}}
\begin{document}
\maketitle
This is the second in a series of articles developing abstract
classification theory for classes that have
a notion of prime models over
independent pairs and over chains.
It deals with the problem of smoothness and establishing the existence
and uniqueness
of a `monster model'.
We work here with a predicate for  a
canonically prime model.  In a
forthcoming paper,
entitled, `{Abstract classes with few models have
`homogeneous-universal' models}',
we show how to drop this predicate from the set of basic notions and
still obtain results analogous to those here.

Experience with both first order logic
and more general cases has shown the advantages of working within a
`monster' model that is both `homogeneous-universal' and `saturated'.
Fra\"{i}ss\'e \cite{Fraisse}
for the countable case and
J\'onsson \cite{Jonssonhomuniv}
for arbitrary cardinalities gave algebraic conditions on a
class {\bf K} of models that guaranteed the existence of a model that is
homogeneous and universal for {\bf K}.
Morley and Vaught \cite{MorleyVaught} showed that if {\bf K} is the
class of models of a first order theory then
the algebraic conditions of homogeneity and universality
are equivalent to model theoretic conditions of saturation.  First order
stability theory works within the fiction of a monster model $\Mscr$.
Such a
fiction can be justified as `a saturated model in an inaccessible
cardinal', by speaking of `class' models,
or by asserting the existence of a function $f$ from cardinals
to cardinals such that
any set of data (collection of models, sets, types
etc.) of
cardinality $\mu$ can be taken to exist in a sufficiently saturated
model of cardinality $f(\mu)$.
This paper culminates by establishing the existence and uniqueness
of such a monster
model for classes {\bf K} of the sort discussed in
\cite{BaldwinShelahprimali} that do not have the maximal number of
models.   We avoid some of the cardinality complications of
\cite{MorleyVaught} by specifying closure properties on the class
of models.

One obstacle to the construction that motivates an important closure
condition
is the failure of `smoothness'.
Is there a unique compatibility class of models embedding a given
increasing chain?
%
It is easy for a model to be
compatibility prime (i.e. prime among all models in a joint
embeddability class over a chain) without being absolutely prime over
the chain.  This `failure of smoothness' is a major obstacle
to the uniqueness of a monster model.
Our principal result here shows that this situation implies
the existence of many models in certain cardinalities.
We can improve the `certain'
by the addition of appropriate set theoretic hypotheses.

The smoothness problem also arose in \cite{Shelahuniversal}.
Even when the
union of a chain is in {\bf K} it does not follow that it can be
{\bf K}-embedded in every member of {\bf K} which contains the chain.
The argument showing this situation implies many models is generalized
here.  However, we have further difficulties.
In the context of \cite{Shelahuniversal}
once full smoothness
(unions of chains are in the class and
are absolutely prime) is established in each cardinality, one
can prove a representation theorem as in \cite{Shelahnonelemii} to recover
a syntactic (omitting types in an infinitary language) definition for
the abstractly given class.  From this one obtains full information
about the Lowenheim-Skolem number of {\bf K} and in particular that
there are models in all sufficiently large powers.  The
examples exhibited in
Section~\ref{smaade}
show this is too much to hope for in our current situation.
Even the simplest case we have in mind, $\aleph_1$-saturated models of
strictly stable theories, gives trouble in $\lambda$ if
$\lambda^{\omega} > \lambda$.
This illustrates one of the added complexities of the more general
situation.  Many properties that in the first order case hold on a final
segment of the cardinals hold only intermittently in the general case.
This greatly complicates arguments by induction and presents the problem
of analyzing the spectrum where a given property holds.


This paper depends heavily on the notations established in Chapters I
and II of \cite{BaldwinShelahprimali}; we do not use the results of
Chapter III.  Reference to \cite{Shelahuniversal} is helpful
since we are generalizing the context of that paper
but we do not
expressly rely on any of the results there.  Some arguments are referred
to analogous proofs in \cite{Shelahuniversal} and
\cite{ShelahMakkaicat}.

Section~\ref{smacpr}
 of this paper recapitulates the properties of canonical prime
models over chains and contains some examples illustrating the efficacy
of the notion.  Sections 2 and 4 fix some basic notations and
assumptions.  Section 2 deals with downward L\"owenheim-Skolem
phenomena; the upwards L\"owenheim-Skolem theorem is considered in
Section 3.  In Section~\ref{smast} we describe the combinatorial
principles used in the paper.  Section~\ref{smadc} introduces a useful
game and some axioms on double chains that are used to show Player I has
a winning strategy if {\bf K} is not smooth.

In Section~\ref{smainv} we consider two important ideas.  First we
discuss the notion `{\bf K} codes stationary sets' --
a particularly strong
form of a nonstructure result for {\bf K}.  Then we consider two
different ways that a model might code a stationary set.  Dealing with
canonically prime models over chains rather than just unions
of chains introduces
subtleties into the decomposition of models that make the process of
`taking points of continuity' more complicated than in earlier studies.

Section~\ref{smasm} contains the main technical results of this paper,
showing that with appropriate set theory, if {\bf K} is
not smooth then it codes stationary sets.

In Section~\ref{smamons} we assume that {\bf K} is smooth.  We are then
able to i) construct and prove the categoricity of
a monster model, ii) introduce the notion of a
type, and iii) recover the Morley-Vaught equivalence of saturation with
homogeneous-universality.
We conclude in Section~\ref{conc} with
a discussion of further problems.
\section{Prime models over chains}
\label{smacpr}

We begin by reviewing the discussion in \cite{BaldwinShelahprimali} of
prime models over chains.

Let $\Mun =\langle M_{\alpha}, f_{\alpha,\beta}
:\alpha <\beta < \delta\rangle$ be an increasing chain of members of
{\bf K}.
An embedding $f$ of $\Mun$ into a structure $M$ is a family of maps
$f_i:M_i \mapsto M$ that commute with the $f_{i,j}$.
As for any diagram, there is an equivalence relation of
`compatibility over $\Mun$'.
Two triples
$(\Mun, f, M)$  and
$(\Mun, g, N)$,
where $f$ ($g$) is a {\bf K}-embedding
of $\Mun$ into $M$ ($N$),  are compatible if there exists an $M'$ and
$f_1$ ($g_1$) mapping $M$ ($N$) into $M'$ such that $f_1 \ocirc f$ and
$g_1 \ocirc g$ agree on $\Mun$ (i.e. on each $M_i$).
This relation is transitive since {\bf K} has the amalgamation property.

Now $M$ is {\em compatibility prime} over $(\Mun,f)$ if it can be
embedded
over $f$ into every model compatible with it.
In Section II.3 of \cite{BaldwinShelahprimali} we introduced the
relation $\cpr$ for canonically prime,
characterized it by an axiom Ch1,  and then
asserted the existence of canonically prime models by axiom Ch2.
Before stating the basic characterization and  existence axioms used in
this paper we need some further notation.

\jtbdef \mbox{}
\label{esscont}
\begin{enumerate}
\item
The chain
$\langle M_i, f_{i,j} :i < j < \beta\rangle $
is  {\em essentially
{\bf K}-continuous
at} $\delta < \beta$
if
there is a model $M'_{\delta}$
that is canonically prime over
$\Mun_{\delta}$ and compatible with
$M_{\delta}$ over
$\langle M_i, f_{i,j} :i < j < \delta\rangle $.

\item
The chain $\Mun$ is  {\em essentially
{\bf K}-continuous}
if for each limit ordinal $\delta < \beta$,  $\Mun$
is   essentially  {\bf K}-continuous
at $\delta$.
\end{enumerate}

We will use the following
slight variants on the existence and characterization axiom in
Section II.3 of \cite{BaldwinShelahprimali}.  Note that by Axiom Ch1',
the $M'_{\delta}$ in Definition~\ref{esscont} can be
$\subm_{\bf K}$-embedded in $M_{\delta}$ over $\Mun|\delta$.

\jtbnumpar{Axioms for Canonically Prime Models}\mbox{}
\begin{description}
\item
[Axiom Ch1'] $\
cpr(\Mun, M,f)$
implies
\begin{enumerate}
\item
$\Mun$ is an essentially {\bf K}-continuous chain,
\item  $M$ is compatibility prime over $\Mun$ via $f$.
\end{enumerate}
\item
[Axiom Ch2']
If $\Mun$ is essentially {\bf K}-continuous there is a
canonically prime model over $\Mun$.
\end{description}

Clearly, if {\bf K} satisfies Ch1' and Ch2' each essentially {\bf
K}-continuous chain can be refined to a {\bf K}-continuous chain.
The following example and Example~\ref{reinonsmooth}
show the necessity of introducing the
predicate $\cpr$
rather than just working with models that satisfy the
definition of compatibility prime.

\jtbnumpar{Example}
\label{omega1tree}
Fix a language with $\omega_1$ unary predicates $L_i$ (for level) and a
binary relation $\prec$.  Let {\bf K} be the collection of structures
isomorphic to structures of the form
$\langle    A, L_i, \prec \rangle$ where
$A$ is a subset of $
^{<\omega_1}
\lambda$  closed under initial segment and containing no uncountable
branch,
$\prec$ is
interpreted as initial segment, and $L_i(f)$ holds if $f\in A$
has length $i$.
Now for $M,N\in {\bf K} $, let $M\subm N$ if $M \subseteq N$ and every
$\omega$-chain in $M$ that is unbounded in $M$ remains unbounded in $N$.

Let $j_{\alpha}$ denote the sequence $\langle k:k<\alpha\rangle$ and
let
$M_i$ be
the member of {\bf K}
whose universe  is $\{j_{\alpha}:\alpha <i\}$.
Then {$\Mun= \langle M_i:i<\omega_1 \ {\rm and } \ i $}
 is
not a limit ordinal $\rangle$ is a {\bf K}-increasing
chain
of members of {\bf K}.
(We do not include the $M_i$ for limit $i$, since if $i$ is a limit
ordinal $M_i$ is not a {\bf K}-submodel of $M_{i+1}$.)
Now in the natural sense for `compatibility prime models' this chain is
continuous.  For each limit $i$, $M_{i+1}$ is compatibility prime in
its compatibility class over $\Mun|(i+1)$ -- the compatibility class of
models with a top on the chain.  But the union of this chain is not a
member of {\bf K} and has no extension in {\bf K}.

So for this example if we tried to introduce `prime'-models over chains
by definition, Axiom Ch2' would fail.
If
in this context
we define $\cpr$
to mean union then the
chain is not even essentially {\bf K}-continuous
and so axiom Ch2' does not
require the existence of a $\cpr$-model over $\Mun$.

Here is an example where the definition of cpr is somewhat more
complicated.

\jtbnumpar{Example} Let {\bf K} be the class of triples
$\langle T,<,Q \rangle$ where $T$ is
a tree
partially ordered by $<$
that has $\omega_1$ levels and
such that each
increasing $\omega$ sequence has a unique least upper bound.  $Q$ is a
unary relation on $T$ such that if $\{t(i):i < \omega_1\}$ is an
enumeration in the tree order of a branch then $\{i:t(i) \in Q\}$ is not
stationary.  For $M, N \in {\bf K} $, write $M \subm N$ if $M$ is a
substructure of $N$ in the usual sense and each element of $M$ has the
same level (height) in $N$.

If $\Mun$ is an increasing chain of {\bf K}-models of length $\mu$, the
canonically prime model over $\Mun$ will be the union of the chain if
$\cf(\mu)$ is uncountable.  If $\cf (\mu) = \aleph_0$, the canonically
prime model will be the union plus the addition of limit points for
increasing $\omega$-sequences but with no new elements added to $Q$.

These examples may seem sterile.  Note however, one of the achievements
of first order stability theory is to reduce the structure of quite
complicated models to trees $\lambda^{<\omega}$.  It is natural to
expect that trees of greater height will arise in investigating
infinitary logics.  Moreover, it is essential to understand these
`barebones' examples before one can expect to deal with more complicated
matters.
In particular, in this framework we expect to discuss the class
of $\aleph_1$-saturated models of a strictly stable theory.
We can not expect to reduce the structure of models of such theories
to anything simpler than a tree with $\omega_1$-levels.
The following example provides another reason for introducing the
predicate $\cpr$.

\jtbnumpar{Example}
\label{reinonsmooth}
Let $T$ be the theory {\bf REI$_{\omega}$} of countably many refining
equivalence relations with infinite splitting \cite[page
81]{Baldwinbook}.
Let {\bf K} be the class of $\aleph_1$-saturated models of $T$ and
define $M \subm N$ if no $E_{\omega}$-class of $M$ is extended in $N$.
($E_{\omega}$ denotes
the intersection of the $E_i$ for finite $i$.)
Now there are many choices for the interpretation of the predicate
$\cpr$, namely the $\kappa$-saturated prime model for each uncountable
$\kappa$.  (Note that if $\Mun = \langle M_i:i < \omega \rangle$  the
models prime among the
$\kappa$ and $\mu$ saturated models respectively containing $\cup \Mun$
are incompatible over $\Mun$ if $\mu \neq \kappa$.)

Thus, the canonically prime model becomes canonical only with the addition of t
he predicate cpr.  There are a number of reasonable candidates in the basic lan
guage and we have to add a predicate to distinguish one of them.
The last
example shows that we should demand that
cpr models are compatible.  This property was not needed in
\cite{BaldwinShelahprimali} but we need it here to prove smoothness.
Its significance is explained in Paragraph~\ref{semismooth}.

\jtbnumpar{Axiom Ch4}
Let $\Mun$ be a {\bf K}-continuous chain and suppose both $\cpr(\Mun,M)$
 and $\cpr(\Mun,N)$ hold.  Then $M$ and $N$ are compatible over $\Mun$.
\section{Adequate Classes}
\label{smaade}

This paper can be considered as a reflection on the construction
of a homogeneous universal model as in
Fraisse \cite{Fraisse}, J\'onsson \cite{Jonssonhomuniv}, and Morley and
Vaught \cite{MorleyVaught}.  These constructions
begin with a class {\bf K} that
satisfies the amalgamation and joint embedding properties.
They
have assumptions of two further sorts: L\"{o}wenheim-Skolem properties
and closure under unions of chains.

We deal with these assumptions in
two ways.  Some are properties of the kinds of classes
we intend to study; we just posit them.
For others we are able to establish within our
context a dichotomy between the property holding and
a nonstructure result for the class.  Most of this paper
is dedicated to the second half of the dichotomy; in this
section we sum up the basic properties we are willing to
assume.

We begin by fixing the language.

\jtbnumpar{Vocabulary}  Recall that each class {\bf K} is
a
collection of structures of fixed vocabulary (i.e.similarity type)
$\tau_{{\bf K} }$. We define
a number of invariants below.  We will require that the
cardinality of
$\tau_{{\bf K} }$ is less than or equal to any of our invariants.
If we did not make this simplifying
assumption we would have to modify each invariant to the maximum of
the current definition and
$|\tau_{{\bf K} }|$.
This would complicate the notation but not affect the arguments in any
essential way.

As usual we denote by {\bf K}$_{\lambda}$ ({\bf K}$_{<\lambda}$ ) the
collection of members of {\bf K} with cardinality $\lambda$
($<\lambda$).
In the next axiom we introduce a cardinal $\chin$.

\jtbnumpar{$\chin$ introduced}
\mbox{}
\begin{description}
    \item[Axiom S0]
$\chin$ is a regular cardinal greater than or equal to $|\tau_{{\bf
K}}|$. \end{description}

Now let us consider L\"{o}wenheim-Skolem phenomena.  In the first
order case,
the upwards L\"{o}wenheim-Skolem property is deduced from the
compactness theorem;
the downwards L\"{o}wenheim-Skolem property holds by the
ability to form elementary submodels by adding {\em finitary} Skolem
functions.
In Section~\ref{smalowsk} we show that an upwards
L\"{o}wenheim-Skolem property can be derived from the basic
assumptions of \cite{BaldwinShelahprimali}.

The finitary nature of the Skolem functions in the first order case
guarantees that the hull of a set of power $\lambda> \chin$
has power
$\lambda$.
Since we now deal
with  essentially infinitary functions, we cannot
make this demand for all $\lambda$.
If there are $\kappa$-ary functions it is likely to fail in cardinals of
cofinality $\kappa$.
We assume a downwards
L\"{o}wenheim-Skolem property in many but not all cardinals.
We justify this assumption in
two ways.  First the condition holds for the classes (most importantly,
$\aleph_1$-saturated models of strictly stable theories) that we intend
to consider.  Secondly, the assumption holds for any class where the
models can be generated by $\kappa$-ary Skolem functions for some
$\kappa$ that depends only on {\bf K} and the similarity type.

\jtbdef
{\bf K} has the $\lambda$-L\"{o}wenheim-Skolem property if for
each
$M\in {\bf K} $ and $A\subseteq M$ with $|A| \leq \lambda$ there exists
an $N$ with $A \subseteq N \subm M$ and $|N| \leq \lambda$.

Replacing the two occurences
of $\leq\lambda$ in the definition of the
$\lambda$-L\"{o}wenheim-Skolem property by $<\lambda$ we obtain
the $(<\lambda)$-L\"{o}wenheim-Skolem property.  If $\mu = \lambda^+$
then the $\lambda$-L\"{o}wenheim-Skolem property and the
$(< \mu)$-L\"{o}wenheim-Skolem property are equivalent.

Note that {\bf K} may have the $\lambda$-L\"{o}wenheim-Skolem property
and fail to have the $\lambda'$-L\"{o}wenheim-Skolem property for some
$\lambda' > \lambda$.

LS({\bf K}) denotes the least $\lambda$ such that {\bf K} has the
$\lambda$-L\"{o}wenheim-Skolem property.

\jtbnumpar{Downward L\"{o}wenheim-Skolem property}
\label{dlsp}
\mbox{}
\begin{description}
    \item[Axiom S1]
There exists a $\chi$
such that for every $\lambda$, if
$\lambda^\chi =\lambda$ then {\bf K} has the $\lambda$-L\"{o}wenheim
Skolem property.
    \item{Notation}
$\chin$ denotes
the least such $\chi$.
$\chik = (\sup(\chin,LS({\bf K}))^+$.
\end{description}

\jtbnumpar{Example}
\label{reinonsmooth1}
Examination of Example~\ref{reinonsmooth} shows
that as stated it does not satisfy the $\lambda$-L\"owenheim-Skolem
property for
any $\lambda$.  An appropriate modification is to consider the class
{\bf K}$^{\mu}$ of models of $T$ that are $\aleph_1$-saturated but each
$E_{\omega}$ class has less than $\mu$-elements.  Then {\bf K}$^{\mu}$
satisfies the $\lambda$-L\"{o}wenheim-Skolem property for any $\lambda
\geq \mu$
and we are able to
apply our main results.

We easily deduce from the $<\lambda$-L\"{o}wenheim-Skolem property the
following decomposition of members of {\bf K}$_\lambda$.  Note that no
continuity requirement is imposed on the chain.

\begin{prop}
\label{decomp1}
If
{\bf K} satisfies the $<\lambda$-L\"{o}wenheim-Skolem
property,
$\lambda$ is regular,
and $M \in {\bf K} $ has cardinality $\lambda$ then $M$ can be
written as $
\bigcup_{i < \lambda} M_i$ where each $M_i$ has power less than
$\lambda$ and $M_i \subm M_j \subm M$ for $i<j<\lambda$.
\end{prop}

We describe chains by a pair of cardinals (size, cofinality) bounding
the size of the models in the chain and the cofinality of the
chain.

\jtbnot A $(\lambda,\kappa)$-chain is a {\bf K}-increasing chain
($i<j$ implies $M_i \subm M_j$)
of
cofinality $\kappa$ of {\bf K}-structures which each have cardinality
$\lambda$.
\par
We define in the obvious way variants on notations of this sort such
as a $(<\lambda,\kappa)$ chain.  Unfortunately, different decisions
about $<$ versus $\leq$ are required at different points and the
complications of notation are needed.
\jtbdef
\begin{enumerate}
\item A chain $\Mun$ is bounded if
for some $M\in {\bf K} $ there is a {\bf K}-embedding of
$\Mun$ into $M$.
\item
{\bf K} is {\em
$(\leq\lambda,\kappa)$-bounded} if  every
$(\leq\lambda,\kappa)$-chain is bounded.
\end{enumerate}
\medskip

To assert {\bf K} is
$(\leq\lambda,\kappa)$-bounded
imposes a nontrivial condition even in the
presence of axiom Ch2' because there is no continuity
assumption on the chain.
Moreover, a demand of boundedness is not comparable to a demand for the
L\"{o}wenheim-Skolem property; it is a demand that a certain
abstract diagram have a concrete realization.
It is easy to construct examples of abstract classes where boundedness
fails if there is a bound on the size of the models in the class.
We describe several more interesting examples in
Paragraph~\ref{variousexamples}.

\jtbnumpar{Alternatives}  Here is a natural further notion.
{\bf K} is
{\bf K}-{\em weakly bounded}
(with appropriate parameters)
if every
{\bf K}-continuous
chain is bounded.
We will not actually have to consider this notion because the
existence of canonically prime models implies that {\bf K} is
{\bf K}-weakly
bounded.
\jtbdef  The cardinal $\lambda$ is {\em {\bf K}-inaccessible} if it
satisfies the following two conditions.
\label{definacc}
 \begin{enumerate}
    \item Any free amalagam $\langle M_0,M_1,M_2\rangle$ with $|M_1|,
|M_2| < \lambda$ can be extended to $\langle M_0,M_1,M_2,M_3\rangle$
with $|M_3| < \lambda$.
    \item  Any $(<\lambda,< \lambda)$ chain, which is bounded, is
bounded by a model with cardinality less than $\lambda$.
\end{enumerate}
\endjtbdef

Since $\lambda^+$ is {\bf K}-inaccessible if {\bf K} satisfies the
$\lambda$-L\"{o}wenheim-Skolem property
we deduce
from S1 the following proposition.
The first clause shows there are an abundance of {\bf K}-inaccessible
cardinals.
For many of the results of this paper it suffices for $\lambda$ to
be {\bf K}-inaccessible rather than requiring the
$<\lambda$-L\"{o}wenheim-Skolem property.

\begin{lemm}
\label{existinacc1}
Suppose $\lambda$ is
greater than $\chi_{\bf K}$ and {\bf K} satisfies S1.
\begin{enumerate}
\item
If $\lambda^{\chin} = \lambda$
then $\lambda^+$
is {\bf K}-inaccessible.
\item If
$\lambda$ is a strongly
inaccessible cardinal
then $\lambda$
is {\bf K}-inaccessible.
\end{enumerate}
\end{lemm}

The following examples show that some of
the classes we want to investigate
have models only in intermittent cardinalities.

\jtbnumpar{Examples}
\label{sporadicspectrum}
Let {\bf K} be the class of
$\aleph_1$-saturated models of a countable
strictly
stable theory $T$.
\begin{enumerate}
\item
If $\lambda^{\omega}> \lambda$
then there are sets
with power $\lambda$ which are contained
in no member of {\bf K}
with power $\lambda$.

\item  If  $T$ is the model completion of the theory of countably many
unary functions
there
is no member of {\bf K}
with power $\lambda$
if $\lambda^{\omega}> \lambda$.
\end{enumerate}

We modify our notion of adequate class from Chapter III of
\cite{BaldwinShelahprimali} to incorporate these ideas.

\jtbnumpar{Adequate Class}
\label{adequate}
{\bf
We assume in this paper axiom groups A and C,
Axioms D1 and D2 from group D,
(all from \cite{BaldwinShelahprimali}),
Axioms Ch1', Ch2' and Ch4 from Section \ref{smacpr}, and
Axioms S0 and
S1  from this section.  A class {\bf K} satisfying these conditions
is called adequate.}

One of our major uses of the L\"owenheim-Skolem property is to guarantee
the existence of {\bf K}-inaccessible cardinals as in
Lemma~\ref{existinacc1}.  We now note that this conclusion can be deduced
from very weak model theory and a not terribly strong set theoretic
hypothesis.  We begin by describing the set theoretic hypothesis.

\jtbdef We say {\em $\infty$ is Mahlo} if for any class $C$ of cardinals
that is closed and unbounded in the class of all cardinals there is a
weakly
inaccessible cardinal $\mu$ such that $C\inter \mu$ is an unbounded
subset  of $\mu$.

In fact, the $\mu$ of the definition could be taken as strongly inaccessible
since the strong limit cardinals form a closed unbounded class.

\begin{thm}
\label{existinacc2}
Suppose $\infty$ is Mahlo and that {\bf K} is a class
of $\tau$-stuctures
that is closed under isomorphism, satisfies axiom C1
(existence of free amalgamations of pairs) and is
$(<\infty,<\infty)$-bounded.  Then the class of {\bf K}-inaccessible
cardinals is unbounded.   In fact, it has nonempty intersection with any
closed unbounded class of cardinals.\end{thm}

\Proof.  For any cardinal $\lambda$, let $J_1(\lambda)$ be the least
cardinal such that for any $\langle M_0,M_1,M_2\rangle$ with the
universe of each $M_i$ a subset of $\lambda$ and each pair $\langle f_1,
f_2\rangle$ where $f_i$ is a {\bf K}-embedding of $M_0$ in $M_i$ (for $i
= 1,2$) there is an $M_3$
and a pair of maps $g_1, g_2$ that complete the amalgamation with $|M_3| <
J_1(\lambda)$.  Similarly, let $J_2(\lambda)$ be the least cardinal such
that any $(\leq\lambda,\leq\lambda)$-chain is bounded by a model of
cardinality less  than $J_2(\lambda)$.  Finally, let $J(\lambda)$ be the
maximum of $J_1(\lambda), J_2(\lambda)$.  Now the set $C = \{\lambda:\mu
<\lambda \text{ implies } J(\mu) < \lambda \}$ is closed and unbounded.
Since $\infty$ is Mahlo, there is an inaccessible cardinal $\chi$ with
$C \inter \chi$ unbounded in $\chi$.  But then $\chi$ is {\bf
K}-inaccessible.
It is easy to vary this argument to show there are actually a proper
class of {\bf K}-inaccessibles and indeed that that class is
`stationary'.

\section{Upwards L\"{o}wenheim-Skolem phenomena}
\label{smalowsk}

As the examples in the Section~\ref{sporadicspectrum}
show, it is impossible to get
the full
L\"{o}wenheim-Skolem-Tarski phenomenon --- models in all sufficiently
large cardinals --- in the most general situation we are studying.
Nevertheless we can establish an
upwards
L\"{o}wenheim-Skolem theorem.
We show that $\chik$
is a Hanf number for models of
{\bf K}.

These results generalize
(with little change in the proof)
and imply
Fact V.1.2 of
\cite{Shelahuni2}.

\jtbnumpar{Remark}  The real significance of
the following theorem is
that it does not rely on axiom C7 (disjointness).  With that axiom the
second
part of the following theorem is trivial.  We asserted in
\cite{BaldwinShelahprimali} that the use of C7 was primarily to ease
notation; this argument keeps us true to that assertion.

Recall that we have assumed for simplicity that $|\tau_{{\bf K} }| \leq
\chi_{\bf K} $.

\begin{thm}
\label{upls}
Suppose {\bf K} has the $\chi$-L\"{o}wenheim-Skolem
property  and there is a member
$M$ of {\bf K} with cardinality greater than $2^{\chi}$.  Then
\begin{enumerate}
\item There exist
$\langle M_0, M_1, M_2, M_3 \rangle$ such that $NF(M_0,M_1,M_2,M_3)$ and
there is
a nontrivial (i.e. not the identity on $M_1$)
isomorphism of $M_1$ onto $M_2$ over $M_0$.
\item There exist arbitrarily large members of {\bf K}.
\end{enumerate}
\end{thm}

\proof.  The proof of conclusion i) is exactly as in \cite{Shelahuni2}.
That is, since
$|M| > 2^{\chi}$, by the $\chi$-L\"{o}wenheim-Skolem property, we
can fix $M_0 \subm M$ with $|M_0 |\leq \chi$ and
choose for each $c\in M- M_0$ an $N_c$ with $M_0 \subm N_c \subm M$, $c
\in N_c$, and $|N_c| \leq \chi$.  Expand the language $L$ of {\bf K} to
$L'$ by adding names for $\{
d:d \in M_0\}$ and let $L''$ contain one more constant symbol.
There are at most $2^\chi$ isomorphism
types of $L''$-structures
$\langle N_c,c\rangle$
satisfying the diagram of
$M_0$
so there are $c_1\neq c_2\in M$ with
$\langle N_{c_1},c_1\rangle
\iso\langle N_{c_2},c_2\rangle$.  Thus there is an isomorphism
$f$ from $N_{c_1}$ onto $N_{c_2}$ over $M_0$ with $f(c_1) = c_2$.
Applying Axiom C2
(existence of free amalgams), we can choose
$M_3$ and
$g:N_{c_1}\iso M_1 $ over $M_0$
with
$NF(M_0,M_1,M,M_3)$.
Now by monotonicity we have both
$NF(M_0,M_1,N_{c_1},M_3)$
and
$NF(M_0,M_1,N_{c_2},M_3)$.  Let $c$ denote $g(c_1)$.  Now not both
$g^{-1}(c)= c_1$
and $f\circ g^{-1}(c)=c_2$ equal $c$.  So one of
$N_{c_1}$ and $N_{c_2}$ can serve as the required $M_2$.

Our proof of the existence of arbitrarily large models actually only
relies on conclusion~i).  Let $c \in M_1$ be such that the isomorphism
$f$ of $M_1$ and $M_2$ moves $c$.  For any $\lambda$, we define by
induction on $\alpha \leq \lambda$ a {\bf K}-continuous sequence of
models $M^{\alpha}$ such that $|M^{\lambda}|\geq \lambda$ as required.
As an auxiliary in the construction we define $f^{\alpha}$ and
$N^{\alpha}$ such that $f^{\alpha}$ is a nontrivial isomorphism between
$M_3$ and $N^{\alpha}$.  We demand
$NF(M_0,N^{\alpha},M^{\alpha},M^{\alpha +1})$.

For $\alpha = 0$,
let $M^0 = M_3$.  At stage $\alpha + 1$ we define $f^{\alpha}$,
$N^{\alpha}$ and
$M^{\alpha+1}$
by invoking the existence axiom to obtain:
$NF(M_0,N^{\alpha},M^{\alpha},M^{\alpha +1})$ and $f^{\alpha}:M_3
\mapsto N^{\alpha}$.
For limit $\alpha$, choose $M^{\alpha}$ canonically prime over its
predecessors.

To obtain the cardinality requirement it suffices to show that if
$\alpha< \lambda$ then $f^{\alpha}(c) \not\in M^{\alpha}$.
Fix $\alpha$ and let
$A_1$ denote $f^{\alpha}(M_1)$
and
$A_2$ denote $f^{\alpha}(M_2)$.
We have
$NF(M_0,A_1,M^{\alpha},M^{\alpha +1})$
and
$NF(M_0,A_2,M^{\alpha},M^{\alpha +1})$
by  the construction.
Again from the construction $g^{\alpha} = f^{\alpha}\restrict M_2 \circ
f \circ (f^{\alpha}\restrict M_1)^{-1}$ is an isomorphism between $A_1$
and $A_2$ over $M_0$.  By the weak uniqueness axiom (C5) (see Lemma
I.1.7 of \cite{BaldwinShelahprimali}), $g^{\alpha}$
extends to an isomorphism $g_{\alpha}$ between $A_1$ and $A_2$ which
fixes $M^{\alpha}$ pointwise.  Now, if $f^{\alpha}(c) \in M^{\alpha}$,
$g^{\alpha}(f^{\alpha}(c))
=
g_{\alpha}
(f^{\alpha}(c))
=
(f^{\alpha}(c))$.
But
$g^{\alpha}(f^{\alpha}(c))
= f^{\alpha}(f(c))$ (by the definition of $g^{\alpha}$) so $f$ fixes
$c$.  This contradiction yields conclusion ii).
\endproof
Noticing that
the existence of a nontrivial map implies the existence of
a nontrivial amalgamation and that only conclusion~i) was used in the
proof of conclusion~ii) we can reformulate the theorem as follows.

\begin{cor}  If
{\bf K} does not have arbitrarily large models then all
members $M$ of {\bf K} have cardinality less than $\chik$.
Moreover, if $N \subm M \in {\bf K} $, there is no nontrivial
automorphism of $M$ fixing $N$.
\end{cor}

\Proof.  Note the definition of $\chik$ (\ref{dlsp})
and apply Lemma~\ref{upls} with
$\chi$ as $\ls(K)$.
Thus
the models of a class with a bound on the size of its models are
all `almost rigid'.
These arguments give some more local information.

\jtbdef
The structure $M$ is a {\em maximal} model in {\bf K} if there is no
proper {\bf K}-extension of $M$.

\begin{cor}
\label{nomax}
If $|M| > 2^{\chi}$ and {\bf K} has the $\chi$-L\"owenheim-Skolem
property then $M$ is not a maximal model in {\bf K}.
Thus if $|M|\geq \chik$, $M$ is not a maximal model.
 \end{cor}

\section{Tops for chains}
\label{smatops}

We discuss in this section several requirements on a
model
that bounds a chain.
Shelah has emphasized (e.g. \cite{Shelahnonelemii,Shelahuniversal})
that the Tarski
union theorem has two aspects.  One is the assertion that the union of
an elementary chain is an elementary extension of each member of the
chain and thus a member of any elementary class containing the chain;
the second is the assertion that the union is an elementary
submodel of any elementary extension of each member of the chain.
First we consider the second aspect.

\jtbdef
\begin{enumerate}
\item
The class {\bf K} is
{\em $(<\lambda,\kappa)$-smooth}
if there is a unique compatibility class over every
$(<\lambda,\kappa)$-chain.
\item
{\bf K} is {\em smooth} if it $(<\infty,<\infty)$-smooth.
\end{enumerate}

Note that if
the class {\bf K} is
{\em $(<\lambda,\kappa)$-smooth}
then the canonically prime model over
any essentially {\bf K}-continuous $(<\lambda,\kappa)$ chain is
absolutely prime.
Moreover, if {\bf K} is smooth every {\bf K}-increasing chain is
essentially {\bf K}-continuous.

The next two notions
represent the first aspect of
Tarski's theorem; the third unites both aspects.

\jtbnumpar{Closure under unions of chains}
\label{defclosed}
\mbox{}
\begin{enumerate}
    \item
{\bf K} is
{\em $(<\lambda,\kappa)$-closed}
if for any
$(<\lambda,\kappa)$-chain
$\Mun$ inside $N$
the union of $\Mun$ is in {\bf K}
and for each $i$, $M_i \subm \cup \Mun$.
    \item
{\bf K} is
$(<\lambda,\kappa)$-{\em weakly closed}
if for any
$(<\lambda,\kappa)$ {\bf K}-continuous chain
$\Mun$ inside $N$
the union of $\Mun$ is in {\bf K}
and for each $i$, $M_i \subm \cup \Mun$.
\item
{\bf K}
is  {\em fully}
$(<\lambda,\kappa)$-{\em smooth} if the union of every
{\bf K}-continuous
$(<\lambda,\kappa)$-chain  inside $M$ is in {\bf K} and
is absolutely prime over the chain.
\end{enumerate}

The `inside $N$' in i) is perhaps misleading.  We have not asserted
$N\in {\bf K}$, so this is not an a priori assumption of boundedness.
In fact $N$ must exist as the union of $\Mun$.  If {\bf K} is
$(\lambda,\kappa)$-closed it is
$(\lambda,\kappa)$-bounded as the union serves as the bound.

If {\bf K} is
{\em $(<\lambda,\kappa)$-smooth}
and
{\em $(<\lambda,\kappa)$-weakly closed}
then for any
$(<\lambda,\kappa)$ {\bf K}-continuous chain
$\Mun$ inside $N$
the union of $\Mun$ is the canonically prime model over $\Mun$
and
{\bf K}
is  {\em fully}
$(<\lambda,\kappa)$-{\em smooth}.

If {\bf K} is the class of $\aleph_1$-saturated models of a strictly
stable countable theory, {\bf K} is not closed under unions of countable
cofinality but is closed under unions of larger cofinality.
This property of a class being closed under unions of chains with
sufficiently long cofinality is rather common.  For example any class
definable by Skolem functions with infinite but bounded arity will have
this property.
We can rephrase several  properties of some of the examples in
Section~\ref{smacpr} in these terms.

\jtbnumpar{Examples}
\label{variousexamples}
\begin{enumerate}
\item
The class
{\bf K}  of
Example~\ref{omega1tree}
is
$(\infty,\geq \aleph_2)$-bounded,
$(<\infty,\geq \aleph_2)$-closed,
and even
fully $(<\infty,\geq \aleph_2)$-smooth.
But {\bf K} is not $(<\aleph_1,\aleph_1)$ bounded and
not $(<\aleph_1,\aleph_0)$
or $(<\aleph_1,\aleph_1)$-smooth.

\item
Example~\ref{reinonsmooth} shows that for {\bf K} the class of
$\aleph_1$-saturated models of $\REI_{\omega}$ and a particular choice
of $\subm$, the class {\bf K} is not smooth.
For, e.g.,
the prime $\aleph_1$-saturated model over a chain
and the prime $\aleph_2$-saturated model over the same chain
may be incompatible.

Note that for any strictly
stable countable theory and any uncountable $\kappa$, if {\bf K} is the
class of $\kappa$-saturated models of a countable
strictly stable theory, $\subm$
denotes elementary submodel, and $\cpr$ means prime among the
$\kappa$-saturated models then {\bf K} is smooth.  In this case
$\chi_1$ is $\aleph_0$.

\item
Consider again the class
{\bf K}$^{\mu}$
discussed in
Example~\ref{reinonsmooth1}.
The union of a countable {\bf K}-chain may determine an $E_{\omega}$
class that is not realized.  So
{\bf K}$^{\mu}$
is not $(<\mu,\omega)$-closed.  But it is
{$(<\mu,[\aleph_1,<\mu])$-closed} and
$(<\infty,\infty)$-bounded.
\end{enumerate}

Our basic argument will establish a dichotomy between
the following weakening of smoothness and a nonstructure theorem.

\jtbnumpar{Semi-smoothness}
\label{semismooth}
The class {\bf K} is
$(<\lambda,\kappa)$-{\em semismooth} if
for each
$(<\lambda,\kappa)$-{\bf K}-continuous chain $\Mun$
each compatibility class over $\Mun$ contains a canonically prime
model over $\Mun$.

The distinction between smooth and semi-smooth
quickly disappears in the presence of Axiom Ch4.

\begin{lemm}  If {\bf K} is semismooth and satisfies Ch4 then {\bf K} is
smooth.
\end{lemm}
\label{semismoothproof}
\proof.  Ch4 asserts that all $M$ satisfying $\cpr(\Mun,M)$ are
compatible.  Since each compatibility class contains such an $M$ there
is only one compatibility class.

\jtbnumpar{Remark}  By an argument similar to the main results of this
paper (but much simpler) we can show
for a proper class
of $\lambda$
that a class {\bf K} that has prime models over independent pairs
and is closed under unions of chains (of any length)
is fully $(<\lambda,<\chin)$-smooth
unless
{\bf K} codes stationary subsets of $\lambda$ (See Section~\ref{smainv}.)
To establish this result we need the axioms about independence of pairs
enumerated in \cite{BaldwinShelahprimali} and that there is a proper
class of {\bf K}-inaccessible cardinals.  The last condition can be
guaranteed by assuming a L\"owenheim-Skolem property like  Axiom S1
or by assuming $\infty$ is Mahlo as in Theorem~\ref{existinacc2}.

This situation is `half way' between the situation in
\cite{Shelahuniversal} and that considered here.  We replace `closed
under substructure' by the existence of `prime models over independent
pairs' but retain taking limits by unions.

\section{Some variants on $\Box$}
\label{smast}
We discuss in this section some variants
on
Jensen's combinatorial principle $\Box$ which will be useful in model
theoretic applications.
We begin by establishing some notation.

\jtbnot
\begin{enumerate}
\item  For any set of ordinals $C$, $\acc[C]$ denotes the set of
accumulation points of $C$ -- the $\delta \in C$ with $\delta = \sup C
\inter \delta$.  The nonaccumulation points of $C$, $C- \acc[C]$,
are denoted $\nacc[C]$.
\item  For any set of ordinals $S$, $C^{\kappa}(S)$ denotes the elements
of $S$ with cofinality $\kappa$.
\item  For any cardinal $\lambda$, $\sing (\lambda)$, the set of {\em
singular ordinals} less than $\lambda$ is the collection of
limit ordinals less than $\lambda$ that are not regular cardinals.
\item
In the
following $\delta$ always denotes a limit ordinal.
\end{enumerate}

The
following definition is a version to allow singular cardinals
of the $\Box$ principle for $\kappa^+$ in \cite{Devlin}.
We refer to it as `full $\Box$'.
This principle has been deduced only from
strong
extensions of ZFC such as $V =L$.
\cite{Devlin}.
We
will also use here weaker versions, obtained by relativizing to a
stationary set, that are provable in ZFC.

When $\lambda = \mu^+$ Jensen
called the condition here a $\Box$ on $\mu$.  The
version here also applies to limit ordinals and since we will deal with
inaccessibles seems preferable.

\jtbdef
The sequence $\langle C_{\delta}:\delta \in \sing(\lambda)\rangle$
witnesses that $\lambda$ satisfies $\Box$ if satisfies the following
conditions.
\begin{itemize}
\item Each $C_{\delta}$ is a closed unbounded subset of $\delta$.
\item $\otp (C_{\delta}) < \delta$.
\item If $\alpha\in \acc[C_{\delta}]$ then
$C_{\alpha}=C_{\delta} \inter \alpha$.
\end{itemize}
\endjtbdef

Now we proceed to the relativized versions of $\Box$.
The relativization is with respect to two subsets, $S$ and $S^+$.
It is in allowing the $C_{\alpha}$ to be indexed by $S^+$ rather than
all of $\lambda$ that this principal weakens those of Jensen and
can be established in ZFC.

We will consider two relativizations.  In Section~\ref{smasm} we will
see that the two set theoretic principles will allow us two different
model theoretic hypotheses for the main result.  They, in fact
correspond to two different  ways of assigning invariants to models.

\jtbdef
\label{boxa}
We say that $S^+ \subseteq \lambda$ and
$\langle
C_{\alpha}:\alpha \in S^+\rangle$ witness that
the subset $S$ of $C^{\kappa}(\lambda)$ satisfies
$\Box^a_{\lambda,\kappa}(S)$
if $S\subseteq S^+$ and
\begin{enumerate}
\item  $S$ is stationary in $\lambda$.
\item For each $\alpha\in S^+$, $C_{\alpha} \subseteq S^+- S$.
\item  If $\alpha \in S^+$ is not a limit ordinal, $C_{\alpha}$ is a
closed subset of $\alpha$.
\item  If $\delta \in S^+$ is a limit ordinal then
\begin{enumerate}
\item $C_{\delta}$ is a club
in $\delta$.
 \item
$\otp (C_{\delta}) \leq \kappa$.
\item $\otp (C_{\delta})= \kappa$ if and only if $\delta \in S$.
\item All nonaccumulation points of $C_{\delta} $ are successor
ordinals.
\end{enumerate}

\item For all $\beta \in S^+$, if $\alpha \in C_{\beta}$ then
$C_{\alpha} = C_{\beta} \inter \alpha$.
\end{enumerate}

\jtbdef
$\Box^a_{\lambda,\kappa}$
holds if for some subset $S \subseteq \lambda$,
$\Box^a_{\lambda,\kappa}(S)$
holds.

\jtbnumpar{Fact}
\label{justificationboxa}
Suppose $\lambda > \kappa$ are regular cardinals.
If $\lambda$ is a successor of a regular cardinal greater than $\kappa$
or
$\lambda = \mu^+$ and
$ \mu^{\kappa}=\mu$ then
for any stationary
$S \subseteq C^{\kappa}(\lambda)$ there is a stationary $S'\subseteq S$
such that
$\Box^a_{\lambda,\kappa}(S')$ holds.

The proof with $\mu^{\kappa} = \mu$ is on page 276 of
\cite{Shelahsquares}
(see also the appendix to \cite{ShelahMakkaicat});
for regular $\mu$ see Theorem 4.1
in \cite{Shelah351}.


Fact~\ref{justificationboxa} is proved in ZFC; if we want to make
stronger demands on the stationary set $S$ we must extend the set
theory.

\jtbdef The subset $S$ of a cardinal $\lambda$ is said to reflect
in $\delta \in S$ if $S \inter \delta$ is stationary in $\delta$.
We say $S$ reflects if $S$ reflects in some $\delta\in S$.
\endjtbdef

Thus a stationary set that does not reflect is extremely sparse in that
its various initial segments are not stationary.
\jtbnumpar {Fact}
\label{justificationstrongboxa}
If $\lambda > \kappa$ are regular, $\lambda$ is not weakly compact,
and V = L then
for any stationary $S_0 \subseteq C^{\kappa}(\lambda)$ there is a
stationary
$S\subseteq S_0$
that does not reflect
such that
$\Box^a_{\lambda,\kappa}(S)$ holds.

This is a technical variant on the result of \cite{BellerLitman}.
Although this result follows from V=L it is also consistent
with various large
cardinal hypotheses.

We now consider the other relativization of $\Box$.  For it we need a
new filter on the subsets of $\lambda$.

\jtbdef  Let the stationary subset $S$ of the regular cardinal
$\lambda$ index the family of sets $C^* =\{C_{\delta}:\delta \in S\}$.
Then $\id(C^*)$ denotes the collection of subsets $B$ of $\lambda$ such
that there is a cub $C$ of $\lambda$ satisfying: for every $\delta \in
B$, $C_{\delta}$ is not contained in $C$.
We denote the dual filter to $\id(C^*)$ by $\fil(C^*)$.

It is easy to verify that $\id(C^*)$ is an ideal.  Note that $B\not\in
\id(C^*)$ if and only if for every club $C$, there is a $\delta\in B$
with $C_{\delta} \subseteq C$.

\jtbdef
\label{boxb}
[$\Box^b_{\lambda,\kappa,\theta,R}(S,S_1,S_2)$]:
Suppose
$\theta \neq
\kappa, \lambda$ and $R$
are four regular cardinals
with $\theta < \lambda, \kappa^+ < \lambda, R\leq \lambda$ and
$S$ is a subset
of $\lambda$ containing all limit ordinals of cofinality $<R$.
We say
$\Box^b_{\lambda,\kappa,\theta,R}(S,S_1,S_2)$ holds if the following
conditions are satisfied for some
$C^* = \langle C_\delta:\delta \in S\rangle$.
\begin{enumerate}
  \item
$C^* = \langle C_\delta:\delta \in S\rangle$
is a sequence of
  subsets of $\lambda$ satisfying
     \begin{enumerate}
        \item $C_{\delta}$ is a closed subset of $\delta$.
        \item $C_{\delta}\subseteq S$.
        \item If $\delta$ is a limit ordinal then
$C_{\delta}$ is unbounded
        in $\delta$.
        \item If $\delta'$ is an accumulation point of
        $C_{\delta}$ then $C_{\delta'} = C_{\delta}\inter \delta'$.
        \item If $\alpha < \delta_1,\delta_2$ and $\alpha \in
        C_{\delta_1}\inter C_{\delta_2}$ then
        $C_{\delta_1} \inter \alpha  =C_{\delta_2}
        \inter \alpha$.
     \end{enumerate}
  \item  $S_1$ and $S_2$ are disjoint subsets of $C^\kappa(S)$ with
  union $S_0\subseteq S$.
     \begin{enumerate}
        \item If $\delta \in S_0$, $\otp(C_{\delta}) = \kappa$.
        \item If $\beta \in \nacc[C_{\delta}]$ and $\delta \in S_0$
         then $\cf(\beta) = \theta$.
     \end{enumerate}
   \item The ideal $\id(C^*)$ is nontrivial; $S_1$ and
$S_2$ are not in $\id(C^*)$.
 \end{enumerate}

\jtbnumpar{Remark}
\label{explainsquare}
Note that if $\delta \in S_0$ and $\beta \in C_{\delta}$ then
$\cf(\beta) < \kappa$ so $C_{\delta} \inter S_0 = \emptyset$.

\jtbdef
[$\Box^b_{\lambda,\kappa,\theta,R}$]:
We say
$\Box^b_{\lambda,\kappa,\theta,R}$
holds if there exist subsets $S \subset \lambda$ and $C_{\beta}\subseteq
S$ for
$\beta \in S$ that satisfy the following conditions.
\begin{enumerate}
\item S contains each of a family
$\langle T_i:i <\lambda\rangle$
of
sets; each $T_i \subseteq C^{\kappa}(\lambda)$ and the $T_i$ are
distinct
modulo
 $\id(C^*)$.
\item For each $A \subseteq \lambda$ there exist $S_1, S_2\subseteq S$
such that
    \begin{enumerate}
    \item $S_1 = \cup_{i\in A} T_i$ and $S_2 = \cup_{i \not \in A}
T_i$ and
\item
$\Box^b_{\lambda,\kappa,\theta,R}(S,S_1,S_2)$ holds with
$C^* =\langle C_{\beta}:\beta \in S\rangle$.
\end{enumerate}
\end{enumerate}

Now the set theoretic strength required for these combinatorial
principles can be summarised as follows.

\begin{thm}
\label{justificationboxb}
\mbox{}
\begin{enumerate}
\item If $\lambda$ is a successor of a regular cardinal,
$\theta^+ < \lambda$, and
$\kappa^+ < \lambda$
then
$\Box^b_{\lambda,\kappa,\theta,\aleph_0}$
is provable in ZFC.
\item If $\lambda$ is a successor
cardinal,
$\theta < \lambda$, and
$\kappa < \lambda$
then
$\Box^b_{\lambda,\kappa,\theta,R}$
is provable in ZFC + V=L for any $R\leq \lambda$.
\end{enumerate}
\end{thm}

\Proof.
Case 1 is proved in III.6.4, III.7.8\,F\,(3)
of \cite{Shelahuni2}
and in
\cite{Shelah351}.  For
Case 2 consult III.7.8\,G of \cite{Shelahuni2}.


\jtbnumpar{Alternative set theoretic hypotheses}
There are a number of refinements on
conditions sufficient to establish
$\Box^b_{\lambda,\kappa,\theta,R}$.
\begin{enumerate}
\item If $\lambda$ is a successor of a regular cardinal
 or even just `not
Mahlo',
Theorem~\ref{justificationboxb} ii) can be strengthened by replacing
`V= L', by `there is a square on $\lambda$'.
See \cite[III.7.8\,H]{Shelahuni2}.
\item
In fact the conclusion of
Theorem~\ref{justificationboxb} ii) holds for any $\lambda$
that is not weakly compact. (Similar to \cite{BellerLitman}.)
\item
The existence of
a function F such that
$\Box^b_{\lambda,\kappa,\theta,R}$
holds for any regular $\lambda > F(\lambda+\kappa+R)$
and $\theta < \lambda$
is consistent with ZFC + there is a class of supercompact
cardinals.
\end{enumerate}
\section{Invariants}
\label{smainv}

As a first approximation we say a class {\bf K} has a nonstructure
theorem if for many $\lambda$, {\bf K} has $2^{\lambda}$ models of power
$\lambda$.
But this notion can be refined.  For some classes {\bf K} it is possible
to code stationary subsets of regular $\lambda$ by models of {\bf K} in
a uniform and absolute way while other classes have many models for less
uniform reasons.
The distinction between these cases is discussed in
\cite{Shelah200,Shelahbook2nd}.  In the stronger situation we
say, informally, that {\bf K} {\em codes stationary sets}.  We don't
give a formal general definition of this notion, but the two coding
functions we describe below $\sm, \som$ epitomize the idea.

The basic intention  is to assign to each model a stationary set so that
at least modulo some filter on subsets of $\lambda$ nonisomorphic models
yield distinct sets.  Historically (e.g. \cite{Shelahuniversal})
to assign such an invariant one
writes the model $M$ as an ascending chain of submodels and asks for
which limit ordinals is the chain continuous.  The replacement of
continuity by {\bf K}-continuity in this paper makes this procedure more
difficult.  We can succeed in two ways.  Either we add an additional
axiom about canonically prime models and proceed roughly as before or we
work modulo a different filter.  Both of these solutions are expounded
here.

We deal only with classes {\bf K} that are {\em reasonably absolute}.
That is, the property that a structure $M \in {\bf K} $ should be
preserved between {\bf V} and {\bf L} and between {\bf V} and reasonable
forcing extensions of {\bf V}.  Of course a first order class or a class
in a $L_{\infty, \lambda}$ meets this condition (a reasonable forcing in
this context would preserve the family of sequences of length $<\lambda$
of ordinals $\leq
\lambda$).  But somewhat less
syntactic criteria are also included.  For example, if {\bf K} is the
class of $\aleph_1$-saturated models of a strictly stable theory
membership in {\bf K} is preserved if we do not add countable sets of
ordinals.

Clearly,
{\bf K} codes stationary sets implies
{\bf K} has
$2^{\lambda}$ models of power $\lambda$.
But it is a stronger evidence of nonstructure in several respects.
First, the existence of many models is preserved under any forcing
extension that does not add bounded subsets of $\lambda$
and
does not destroy the stationarity of subsets of $\lambda$.  Secondly,
the
existence of many models on a proper class of cardinals is not such a
strong requirement; for example, a multidimensional
(unbounded in the nomenclature of \cite{Baldwinbook}) theory has
$2^{\aleph_{\alpha}}$ models of power $\aleph_{\alpha}$ whenever
$\aleph_\alpha = \alpha$.  However, the class of models of a
first order theory codes stationary sets
only if $T$ is not superstable or has the dimensional order property or
has the omitting types order property.

\jtbdef  \mbox{}
\label{represent}
\begin{enumerate}
\item
A {representation} of a model $M$ with power $\lambda$ (with $\lambda$
regular)
is an
increasing chain $\Mun = \langle M_i:i<\lambda\rangle$ of {\bf
K}-substructures of $M$
such that each
$M_i$ has cardinality less than $\lambda$ and $\cup \Mun = M$.
\item
The representation is {\em proper} if
$\cup \Mun|\delta \subm M$ implies
$M_\delta = \cup \Mun|\delta$.
\end{enumerate}

We showed in Proposition~\ref{decomp1} that if $\lambda$ is a regular
cardinal greater than $\chin$ and {\bf K} satisfies the
$\lambda$-L\"{o}wenheim-Skolem property  then
each model of power $\lambda$ has a
representation.  We will not however have to invoke the
L\"{o}wenheim-Skolem property in our main argument because we
analyze models that are constructed with a representation.
Using Axiom A3,
it is easy to perturb any given representation into a
proper representation.

We now show how to define  invariant functions in our context.  For
the first version we need an additional axiom.

\jtbnumpar{Axiom Ch5}
For every {\bf K}-continuous chain $\Mun = \langle
M_i,f_{i,j}:i,j<\delta\rangle$ and every unbounded $X\subseteq
\delta$, a {\bf
K}-structure $M$ is canonically prime over $(\Mun,f)$ if and only if $M$
is canonically prime over $(\Mun|X,f|X)$.

We refer to this axiom by saying that {\em canonically  prime models
behave on subsequences}.

\jtbdef
Let $\Mun$ be a representation of $M$.
Let $\sm(\Mun,M)$ denote the set of
limit ordinals $\delta
<\lambda$ such that for some $X$ unbounded in $\delta$ a canonically
prime model $N$ over $\Mun|X$ is a {\bf K}-submodel of $M$.

\begin{lemm}
\label{defrep}
If $\Mun$ and $\Nun$ are representations of $M$ and axiom
Ch5 holds then
$\sm(\Mun,M) = \sm(\Nun, M)$ modulo the closed unbounded
filter on $\lambda$.
\end{lemm}

\Proof.  Since $|M|$ is regular, there
is a cub $C$ on $\lambda$ such that every
$\delta \in C$ is
a limit ordinal and
$\cup \Mun|\delta
=\cup \Nun|\delta$ for $\delta \in C$.  Now we claim that for $\delta
\in C$,
$\delta \in \sm(\Mun,M)$ if and only if $\delta\in\sm(\Nun, M)$.
To see this choose an increasing sequence $L_i$ alternately from $\Mun$
and $\Nun$.  Since $\cpr$ behaves on subsequences
the  canonically prime model over the common
subsequence of $\Lun$ and $\Mun$
is a {\bf K}-submodel of $M$
if and only if the canonically prime model over
$\Lun$ is  and similarly for the common subsequence
of $\Lun$ and $\Nun$.  Thus,
$\delta \in \sm(\Mun,M)$ if and only if $\delta\in\sm(\Nun, M)$.

This lemma justifies the following definition.

\jtbdef
Denote the equivalence class modulo $\cub (\lambda)$
of $\sm(M,\Mun)$ for some (any) representation $\Mun$ of $M$
by $\sm(M)$.
We call $\sm (M)$ the {\em smoothness set} of $M$.

We now will describe a second way to assign invariants to models.  This
approach avoids the reliance on Axiom Ch5 at the cost of complicating
(but not increasing the strength of) the set theory.
Recall from Section \ref{smast} the ideal $\id(C^*)$ assigned to a
family of sets
$C^*= \{C_{\beta}:\beta \in S\}$.
Fix for the following definition and arguments subsets $S, S_1, S_2$
satisfying
$\Box^b_{\lambda,\kappa,\theta,R}(S,S_1,S_2)$.
We define a second
invariant function with $C^*$ as a parameter.  It distinguishes models
modulo $\id(C^*)$.

\jtbdef
\label{defsom}
Fix a subset $S$ of $\lambda$ and
$C^*= \{C_{\beta}:\beta \in S\}$.
Let $\Nun$ be a representation of $M$.
$\som(\Nun,C^*,M)$ is the set of $\delta \in S$ such that
\begin{enumerate}
  \item for every $\gamma \in \nacc[C_{\delta}]$, $N_{\gamma} =
\cup_{\alpha < \gamma}N_{\alpha}$,
  \item $\Nun|C_\delta$ is {\bf K}-continuous,
  \item there is an
$N'_{\delta}$ canonically prime over
$\Nun |
\nacc[C_{\delta}]$
that can be {\bf K}-embedded into $M$ over
$\Nun |
\nacc[C_{\delta}]$.
\end{enumerate}

We are entitled to choose an
$N'_{\delta}$ canonically prime over
$\Nun |
\nacc[C_{\delta}]$
in condition iii) because condition i) guarantees
that
$\Nun |
\nacc[C_{\delta}]$ is {\bf K}-continuous.

\begin{lemm}  If $\Mun$ and $\Nun$ are proper representations of $M$
then
$$\som(\Nun,C^*,M)
=\som(\Mun,C^*,M) \text{  modulo } \fil(C^*).$$
\end{lemm}

\Proof. Let
$X_1$ denote
$\som(\Nun,C^*,M)$
and
$X_2$ denote
$\som(\Mun,C^*,M)$.
Without loss of generality we can assume the universe of $M$ is
$\lambda$.
There is a cub $C$ containing only limit ordinals such that for
$\delta\in C$, $\delta =
\cup_{\alpha < \delta}
 M_{\alpha}=
\cup_{\alpha < \delta}
 N_{\alpha}$.


To show $X_1 = X_2$
mod $\fil(C^*)$,
it suffices to show there is a $Y \in {\rm
{FIL}}(C^*)$ such that $X_1 \inter Y = X_2 \inter Y$.  Let $Y =
\{\delta:C_\delta \subseteq C\}$.

   Suppose $\delta \in Y\inter X_1$.  If $\alpha \in \nacc[C_{\delta}]$
then $\alpha \in Y$ implies $\alpha \in C$ which in turn implies
$\alpha =
\cup_{i<\alpha}N_i =
\cup_{i<\alpha}M_i$.  Now $\delta \in X_1$ implies $N_{\alpha} =
\cup_{i<\alpha}N_i$ and $N_{\alpha}\subm M$ so $
\cup_{i<\alpha}N_i \subm M$ and thus
$\cup_{i<\alpha}M_i \subm M$.  But then by properness
$M_\alpha =\cup_{i<\alpha}M_i$.
Thus,
$M_\alpha =
N_{\alpha}$.
That is,
$\Nun|\nacc[C_{\delta}]
=\Mun|\nacc[C_{\delta}]$.  So $\delta \in X_1$ if and only if $\delta \in
X_2$.

In view of the previous lemma we make the following definition.

\jtbdef For any $M$ in the adequate class {\bf K} and some (any) proper
representation $\Mun$ of $M$,
$\som(C^*,M) = (\som(\Mun,C^*,M)/\fil(C^*))$.
\section{Games, Strategies and Double Chains}
\label{smadc}

We will formulate one of the main model theoretic hypotheses
for the major theorem deriving nonstructure from nonsmoothness in
terms of the existence of winning strategies for a certain game.
In this section we describe this game and show how to derive a winning
strategy for it from the assumption that {\bf K} is not smooth.

\jtbdef  A play of Game 1 $(\lambda,\kappa)$
lasts $\kappa$ moves.  Player I chooses models $L_i$ and Player II
chooses models $P_i$ subject to the following conditions.
At move $\beta$,
\begin{enumerate}
\item Player I chooses a model $L_{\beta}$ in {\bf K} of power less than
$\lambda$ that is a proper {\bf K}-extension of all the structures
$P_{\gamma}$ for $\gamma <\beta$.  If $\beta$ is a limit ordinal
less than $\kappa$, $L_{\beta}$ must be chosen canonically
prime over $\langle
P_{\gamma}:\gamma < \beta\rangle$.

\item Player II chooses a model $P_{\beta}$ in {\bf K} of power less
than $\lambda$ that is a {\bf K}-extension of
$L_{\beta}$.
\end{enumerate}

Any player who is unable to make a legal move loses.
Player I wins the game if there is a model $P_{\kappa} \in {\bf
K}$ that extends each $P_{\beta}$ for $\beta < \kappa$ but
the sequence $\langle P_i:i \leq \kappa\rangle$ is not
essentially {\bf K}-continuous.

In order to establish that nonsmoothness implies a winning strategy for
Player I
we need
to consider certain properties of double chains.  We introduce here some
notation and axioms concerning this kind of diagram.

\jtbdef
\begin{enumerate}
\item $\Mun = \langle \Mun^0, \Mun^1 \rangle = \{\langle M^0_i,
M^1_i\rangle :i < \delta\}$ is a {\em double chain} if each $M^0_i \subm
M^1_i$ and $\Mun^0, \Mun^1$ are {\bf K}-increasing
chains.
We say $\Mun$ is (separately)
({\bf K})-continuous
if each of $\Mun^0$
and $\Mun^1$ is
({\bf K})-continuous.
\item
$\Mun$ is a {\em free} double chain
if for each $i< j<\delta$,
$\ind  {M^0_j} {M^0_i} {M^1_i}$ inside $M^1_{j+1}$.

\item $\Mun = \langle \Mun^0, \Mun^1,
\rangle = \{\langle M^0_i,
M^1_j\rangle :i \leq
 \delta+1,j<\delta\}$
  is a {\bf K}-{\em continuous augmented
 double chain} inside $N$ if
$i< \delta$ implies $M^0_i \subm
M^1_i$,
and $\Mun^0, \Mun^1$ are increasing {\bf K}-continuous
chains inside $N$.

\item An augmented double chain
is {\em free} inside $N$ if
for each $i < \delta$,
$$\ind {M^1_i} {M^0_{i}} {M^0_{\delta + 1}}\  \hbox{\rm inside } N.$$
\end{enumerate}
\endjtbdef
We extend the existence axiom Ch2' for  a prime model over a
chain to assert the compatibility of the prime models guaranteed for
each sequence in a double chain.

\jtbnumpar{Axioms concerning double chains}
\mbox{}
\begin{description}
\item [DC1]
  If $\Mun$ is an essentially {\bf K}-continuous free double chain
and $M_1$ is canonically prime over $\Mun^1$ then there is an
$M_0$ that
is canonically prime over $\Mun^0$ such that $M_0$ and $M_1$ are
compatible over $\Mun^0$.

\item[DC2]  If $\Mun$ is an essentially {\bf K}-continuous
free augmented
double chain of length $\delta$ in $M$ then
there is an $N$ with $M\subm N$ and an
$M^1_{\delta} \subm N$ such that
$$\ind {M^1_{\delta}} {M^0_{\delta}} {M^0_{\delta + 1}}
\ \hbox{\rm inside } N $$

and the chain $\Mun^1 \cup \{M^1_{\delta}\}$ is essentially
{\bf K}-continuous.

\end{description}

We will refer to versions of these axioms for chains of
restricted length; we may denote the variant of the axiom for
chains of length less than $\kappa$ as
DCi($<\kappa$).

Note that it would be strictly stronger in DC2 to assert that
$M^1_{\delta}$ is canonically prime over $\Mun^1$ since under
DC2 as stated the canonically prime model over $\Mun^1$ inside
$M^1_{\delta} $
need not contain
$M^0_{\delta}$.

Since we are going to use these axioms to establish smoothness we
indicate some relationships between the properties.  {\bf K} is
$(<\infty,\kappa)$
smooth
means that every {\bf K}-continuous chain of cofinality $\kappa$ has a
single compatibility class over it -- necessarily there will be a
canonically prime model in that class.  DC1 would hold if there were
many compatibility classes over a chain but each had a canonically
prime model (i.e. {\bf K} is semismooth).
In particular it holds at $\kappa$ if {\bf K} is
$(<\infty,\kappa)$-smooth (sometimes read smooth at $\kappa$).
Thus,
the following lemma is easy.
\begin{lemm}
\label{deducedc1}
If {\bf K} is an adequate class that is
$(<\infty,\kappa)$-smooth then
{\bf K} satisfies DC1 for chains
of cofinality $\kappa$.
\end{lemm}

Now we come to the main result of this section.

\begin{lemm}  \label{mainnotsmoothstrat}
Let {\bf K} be an adequate class that is
{$(<\lambda,\kappa)$-bounded}
and suppose {\bf K} is $(<\lambda,<\kappa)$-smooth
but not $(<\lambda,\kappa)$-smooth.
Suppose further that {\bf K} satisfies DC1 and DC2
and $\lambda > \chik$ is {\bf K}-inaccessible.
Then Player I has a winning
strategy for {Game 1} $(\lambda,\kappa)$.
\end{lemm}

\Proof.
Since {\bf K} is $(<\lambda,<\kappa)$ smooth, we can choose a
counterexample
$\Nun = \langle N_i:i <\kappa\rangle$ that is
essentially
{\bf  K}-continuous.  Then $\Nun$ is bounded by two models
$N_{\kappa}$
and $N'_{\kappa}$ that are incompatible over $\Nun$.  If there
exist
$M_{\kappa}$
and $M'_{\kappa}$
canonically prime over $\Nun$ and embeddible
in
$N_{\kappa}$
and $N'_{\kappa}$ respectively,  Axiom Ch4 requires that
$M_{\kappa}$
and $M'_{\kappa}$ are compatible.  But then so are
$N_{\kappa}$
and $N'_{\kappa}$.  From this contradiction we conclude without loss of
generality
that
each $N_i \subm N_{\kappa}$ but that no
canonically prime model over $\Nun$ can be {\bf K}-embedded into
$N_{\kappa}$.
That is, {\bf K} is not semismooth (\ref{semismooth}).
Now Players I and II will choose models
$\langle L_i:1 \leq i <
\kappa\rangle$ and
$\langle P_i:i <
\kappa\rangle$  for a play of Game 1.

We describe a winning strategy for Player I.  The construction requires
auxiliary models $P'_i$, $N^*_i$, and $L'_i$ and isomorphisms
$\alpha_i:L'_i \mapsto L_i$.  They will satisfy the following
conditions.
\begin{enumerate}
\item $\Pun'$ and $\Lun'$ are essentially {\bf K}-continuous  sequences
and the $\alpha_i$
are an increasing sequence of maps.
\item $\ind {P'_i} {N_i} {N_{\kappa}}$
inside $N^*_{i+1}$.
\item $L'_{i+1}$ is  prime over $P'_i \cup N_{i+1}$
inside $N^*_{i+1}$.
\item $\alpha_i$ is an isomorphism between $L'_i$ and $L_i$
mapping $P'_j$ onto $P_j$ for $j<i$.
\item The $N^*_i$ form an essentially {\bf K}-continuous sequence with
$N_i \subm N^*_i$.
\end{enumerate}

Let $L_0= N_0$.
Each successor stage is easy.  Player II has chosen $P_{i}\in {\bf
K}_{<\lambda}$ to extend $L_{i}$.
For Player I's move, apply Axiom D1 (existence of free
amalgamations) to
first choose
$N^*_{i+1}$ to extend
$N^*_{i}$ and $P'_i$ with $P'_i \iso P_i$
by an isomorphism $\hat \alpha_i$ extending $\alpha_i$
and with
$\ind {P'_i} {N_i}
{N_{\kappa}}$ inside $N^*_{i+1}$ to satisfy condition ii).
Then choose
$L'_{i+1}$
to satisfy iii) by the existence of free amalgamations (Axiom D1).
Finally choose $L_{i+1}$
and $\alpha_{i+1}$ extending $\hat \alpha_i$ to satisfy condition iv).
As $\lambda$ is {\bf K}-inaccessible, $N^*_{i+1}$ and $L_{i+1}$ can
be chosen in {\bf K}$_{<\lambda}$.
At a limit ordinal $\delta< \kappa$,
let $\tilde N_{\delta}$ be canonically
prime over
$\langle N^*_i:i < \delta\rangle$.
Then
$\langle \langle N_i:i \leq \delta\rangle \cup \{N_{\kappa}\}$,
$\langle L'_i:i < \delta\rangle\rangle$
is a free augmented double chain inside $\tilde N_{\delta}$.
(Strictly speaking, this is proved by induction on $\beta < \delta$.
Use the base extension axiom to pass from
$\ind {P'_i} {N_i}
{N_{\kappa}}$
to
$\ind {L'_{i+1}} {N_{i+1}}
{N_{\kappa}}$.)
By DC2 there is an $N^*_{\delta}$ {\bf K}-extending $\tilde N_{\delta}$
and an $L'_{\delta}\subm N^*_{\delta}$ with
$\ind {L'_{\delta}} {N_{\delta}}
{N_{\kappa}}$.
Extend $\langle \alpha_i:i<\delta\rangle$ to map $L'_{\delta}$ to
$L_{\delta}$.

Now we show that this strategy wins for Player I.
Since $\Pun$ and $\Pun'$ are isomorphic, it suffices to show that there
is no $P'_{\kappa}$ with $\Pun' \cup\{P'_{\kappa}\}$ essentially {\bf
K}-continuous.  Suppose for contradiction that such a $P'_{\kappa}$
exists.  Since $\langle \Nun \cup \{N_{\kappa}\}, \Pun'\rangle$ is a
free  double chain  inside $N^*_{\kappa}$,
by DC1 the canonically prime model $N'$
over $\Nun$ can be embedded in $P'_{\kappa}$ inside some extension of
$N^*_{\kappa}$.
But then $N_{\kappa}$ and $N'$ are compatible over $\Nun$
contrary to assumption.
\medskip
\jtbnumpar{Remarks}
\begin{enumerate}
\item
Instead of assuming Axiom Ch4 (part of the definition of adequate) we
could have assumed that {\bf K} was not $(<\lambda,\kappa)$-semismooth.
\item
It is tempting to think that by choosing the minimal length of
a sequence witnessing nonsmoothness, we could apply
Lemma~\ref{deducedc1} and avoid assuming DC1.  However,
DC1 is applied for chains of length  $\kappa$ so this ploy
is ineffective.
\item Why is $L_{i+1}$ a proper extension of $P_i$?  Since $L_{i+1}$ was
chosen as an amalgam of $P'_{i}$ and $N_{\kappa}$, this is immediate if
we assume the disjointness axiom (C7).  To avoid this hypothesis we can
demand that each model in the construction have cardinality $>\chik$ and
so not be maximal (by Lemma~\ref{nomax}).  That is why we assumed
$\lambda >\chik$.
\item
Note that DC1 is used to derive the contradiction at the end of the
proof;  DC2 is used to pass through limit stages of the construction.
Thus in the important case when $\kappa = \omega$ we have
\end{enumerate}

\begin{lemm}
\label{notsmoothstrat}
Let {\bf K} be an adequate class
that is
not {$(<\lambda,\omega)$ smooth.}
Suppose  that {\bf K} satisfies DC1.
Then Player I has a winning
strategy for {Game 1} $(\lambda,\omega)$.
\end{lemm}

The choice of $L_i$ according to the winning strategy of Player I
depends only on the sequence $\langle L_j, P_j\rangle$, for $j<i$
(not, for example on some guess about the future of the game).

In the remainder of this section
section we consider a third axiom DC3 on double chains.
The following axiom bears the same relation to DC2 that C5 bears to C2.

\jtbnumpar{Weak uniqueness for prime models over double chains}
\mbox{}
\begin{description}
   \item[DC3]
Suppose that $\Mun$ and $\Nun$ are essentially {\bf K}-continuous
augmented double
chains that are free in $M$ and $N$ respectively and $f$ is an
isomorphism from $\Mun$ onto $\Nun$.  Suppose also that
$\ind {M^1_{\delta}} {M^0_{\delta}} {M^0_{\delta+1}}$ inside $M$
and
$\ind {N^1_{\delta}} {N^0_{\delta}} {N^0_{\delta+1}}$ inside $N$.
Then there is an $\hat M \in {\bf K} $ and {\bf K}-embeddings $h_0$ of
$M$ and $h_1$ of $N$ into $\hat M$ with $h_1 \circ f
= h_0
\restrict \Mun$.
\end{description}

Just as Lemma I.1.8 of \cite{BaldwinShelahprimali} rephrased the weak
uniqueness axiom for amalgamation over vees we can reformulate DC3 as
follows.

\begin{lemm}
\label{dc3use}
Assume DC2 and DC3.
Suppose that $\Mun$ and $\Nun$ are essentially {\bf K}-continuous
augmented double
chains that are free in $M$ and $N$ respectively and $f$ is an
isomorphism from $\Mun$ onto $\Nun$.  Suppose also that
for some $M^1_{\delta} \subm N$,
$\ind {M^1_{\delta}} {M^0_{\delta}} {M^0_{\delta+1}}$ inside $M$.

Then there exist a model $\hat N$ and an isomorphism $h:M \mapsto \hat
N$
such that
$h \supseteq f$ and
$\ind {h(M^1_{\delta})} {N^0_{\delta}} {N^0_{\delta+1}}$ inside
$\hat N$.
\end{lemm}

\begin{quest}  If DC2 and DC3 hold and {\bf K} is not smooth
does Player I have a winning strategy
for Game 1?
\end{quest}
\section{Nonsmoothness implies many models}
\label{smasm}

We show in this section that if the class {\bf K} is not smooth then
{\bf K} codes stationary sets.  These results involve several tradeoffs
between set theory and model theory.  The main result  is proved in ZFC.
Here there are two versions; one uses $\Box^a$ and requires the
hypothesis that cpr behaves on subsequences.  The second uses $\Box^b$
and replaces `cpr behaves on subsequences' with stronger hypotheses
concerning the closure of {\bf K} under unions of chains.  By working in
L  we can reduce our assumptions on which chains are bounded in both
cases.

\subsection{Invariants modulo the CUB filter}

In this subsection we show if {\bf K} is not smooth then for many
$\lambda$ we can code stationary subsets of $\lambda$ by assigning
invariants in the cub filter by the function $\sm$.
Our general strategy for constructing many models is this.
We build a model $M^W$ for each of a family of
$2^{\lambda}$ stationary subsets $W$ of $S$ that are pairwise distinct
modulo the cub filter.
The key point of the construction is that, modulo $\cub(\lambda)$,
we can recover $W$ from
$M^W$ as
$\lambda -\sm(M^W)$.

We need one more piece of notation.

\jtbnot
Fix a square sequence $\langle C_{\alpha}:\alpha \in S\rangle$.
Suppose Player I has a winning strategy for Game 1$(\lambda,\kappa)$.
In the proof of Theorem~\ref{maintheorem:smooth} and some similar later
results we define a {\bf K}-increasing sequence $\Mun$.  We describe
here what is meant by saying a certain $M_{\alpha}$ is chosen by playing
Player I's strategy on $\Mun|C_{\alpha}$.

Let $\langle c_{\beta}:\beta < \beta_0\rangle$ enumerate $C_{\alpha}$.
We regard $\Mun|C_\alpha$ as two sequences
$\langle \Lun,\Pun\rangle$ by
setting for any ordinal $\delta +n$ with $\delta$ a limit ordinal
and $n<\omega$:
$$L_{\delta + n} \text{ is } M_{c_{\delta + 2n}}.$$
$$P_{\delta + n} \text{ is } M_{c_{\delta + 2n +1}}$$
We say $M_{\beta}$ for $\beta = \alpha$ or $\beta \in C_{\alpha}$
is chosen by Player I's winning
strategy on $M|C_{\alpha}$ if the sequence
$\langle \Lun,\Pun\rangle$
associated with $C_{\alpha}\inter \beta \union \{\beta\}$ is
\begin{enumerate}
\item an initial segment of a play of Game 1 $(\lambda,\kappa)$ and
\item Player I's moves in this game follow his winning strategy.
\end{enumerate}

%

Here is the technical  version of the main result with the parameters
and
reliance on the axioms enunciated in
Section~\ref{smaade} and \ref{smatops}
stated explicitly.

Although the assumption that $\lambda$ is {\bf K}-inaccessible is weaker
than the assumption that {\bf K} satisfies the $\lambda$-L\"{o}wenheim
Skolem property it plays the role of the L\"{o}wenheim-Skolem
property in the following construction.
We assume $\lambda > \chik$ and apply Corollary~\ref{nomax} to
avoid the appearance of maximal models in the construction.
\begin{thm}
\label{maintheorem:smooth}
Fix regular cardinals $\kappa < \lambda$.
Suppose the following conditions hold.
\begin{enumerate}
\item
{\bf K} is an adequate class.
\item   Player I has a winning strategy for Game 1 $(\lambda,\kappa)$.
\item
$\lambda$ is a {\bf K}-inaccessible cardinal, for some
stationary $S\subseteq \lambda$,
$\Box^a_{\lambda,\kappa}(S)$ holds and $\lambda > \chik$,
\item
{\bf K} is $(<\lambda,<\lambda)$-bounded.
\item
 $\cpr$ behaves on
subsequences (Ch5).
\item
 {\bf K} is  $(<\lambda,\lambda)$-closed.
\end{enumerate}
Then for any stationary $W\subseteq S$
there is a model $M^W$ and a representation $\Mun^W$
with
$W\subseteq \lambda -\sm(M^W,\Mun^W)$  and
$S^+-W\subseteq \sm(M^W,\Mun^W)$.
\end{thm}

\proof.
Fix
$S^+$ and
$C^* =\langle C_i:i \in S^+\rangle$ to witness
$\Box^a_{\lambda,\kappa}(S)$.
Without loss of generality, $0 \in C$.
Fix also a stationary subset $W$ of $S$.
For each $\alpha < \lambda$ we define a model $M^W_{\alpha}$.
The model $M^{W} = \bigcup_{\alpha <\lambda} M^W_{\alpha}$
constructed in this way is in {\bf K} by condition vi) and
will satisfy the conclusion.

Each of these conditions  depends indirectly on $W$, but
since we are constructing each $M^W$ separately,
we suppress the dependence on $W$ to avoid notational confusion in the
construction.

For each
$\alpha< \lambda$ we define
$M_{\alpha}$
to satisfy the
following requirements.
\begin{enumerate}
\item $|M_{\alpha}| \geq \chik$ (to avoid maximal models).
\item
$\Mun$ ($= \Mun^W$)
is an increasing sequence of members of $K_{<\lambda}$ which is
essentially
{\bf K}-continuous at $\delta$ if $\delta  \in (S^+-W)$.
\item If  $\alpha\in W$ then $\Mun$ is
not essentially {\bf K}-continuous at $\alpha$.

\item $M ( =M^{W}) = \bigcup_{\alpha < \lambda} M_{\alpha}$.
\end{enumerate}
\par

The construction proceeds by induction.  There are a number of cases
depending on whether $\alpha\in W, S$ etc.

\begin{description}
\item[Case I.]
$\alpha \in (W \cup \bigcup_{\delta\in W} C_{\delta})$:
$M_{\alpha}$ is chosen by Player I's winning strategy for
$\Mun|C_{\beta}$
for any $\beta$ with $\alpha \in C_{\beta}$.  (The choice of $\beta$ does
not matter because of the coherence condition in the definition of
$\Box^a_{\lambda,\kappa}$.)
\item[Case II.]
$\alpha \in S^+- (W \cup \bigcup_{\delta\in W} C_{\delta})$:
Then $C_{\alpha} \subseteq S^+-S$ so $\Mun|C_{\alpha}$
is {\bf K}-continuous.
(Remember that Player I plays canonically prime models at limit ordinals of
cofinality less than $\kappa$.)
Choose
$M_{\alpha}$ to be canonically prime over $\Mun|C_{\alpha}$ (which is
the same as canonically prime over $\Mun|\alpha$ since cpr behaves on
subsequences).
\item[Case III.] $\alpha \not \in S^+$.  Choose $M_{\alpha}$ to
bound $\Mun|\alpha$ by $(<\lambda,<\lambda)$ boundedness and with
$|M_{\alpha}| < \lambda$ since $\lambda$ is {\bf K}-inaccessible.
\item[Case IV.]  Any successor ordinal not already done.
Say, $\beta = \gamma+1$. Choose $M_{\beta}$
as a proper {\bf K}-extension of $M_{\alpha}$ by
Corollary~\ref{nomax}.
\end{description}

The cases in the construction are easily seen to be disjoint
(using ii) of the definition of $\Box^a_{\lambda,\kappa}$)
and
inclusive.  If $\delta\in W$ is a limit ordinal
the canonically prime model over
$\Mun|\nacc[C_{\delta}]$ is not compatible with $M_\delta$
since Player I played a winning strategy on $M|C_{\alpha}$.
So, since cpr
behaves on subsequences, neither is the canonically prime model on
$\Mun|A$ for any $A$ unbounded in $\delta$.  Thus $\delta \not \in
\sm(\Mun,M^W)$.  All other limit $\delta\in S^+
$ are in $\sm(\Mun,M^W)$ and we
finish.
\medskip

The next theorem rephrases Theorem~\ref{maintheorem:smooth}
to avoid technicalities.  It shows
that reasonable {\bf K} that are not smooth have many models in all
sufficiently large successor cardinals.  In fact we have the stronger
result that {\bf K} codes stationary subsets of such cardinals.
\begin{thm}[ZFC]
\label{firstquotable}
Let
{\bf K} be an
adequate class and
suppose that {\bf K} satisfies DC1, DC2 and
cpr behaves on subsequences (Ch5).
Suppose there exist $\kappa, \lambda_1$ with
$\kappa< \lambda_1$ such that
{\bf K} is not $(\lambda_1, \kappa)$-smooth.
Then for every {\bf K}-inaccessible $\lambda > \sup(\chik,\lambda_1)$
such that
\begin{enumerate}
\item $\lambda$ is a successor of a regular cardinal,
\item {\bf K} is $(<\lambda,<\lambda)$-bounded,
\item {\bf K} is $(<\lambda,\lambda)$-closed,
\end{enumerate}
{\bf K} has $2^{\lambda}$ models in power $\lambda$.
\end{thm}

\proof.
{\bf K} is not $(<\lambda_{1},\kappa)$-smooth trivially implies
{\bf K} is not
$(<\lambda,\kappa)$-smooth.
We assumed DC1 and DC2 so
by Lemma~\ref{mainnotsmoothstrat}, Player I
has a winning strategy in Game 1 $(\lambda,\kappa)$.
By Fact~\ref{justificationboxa},
there is a stationary $S\subseteq C^{\kappa}(\lambda)$ such that
$\Box^a_{\lambda,\kappa}$ holds.  The result now follows from the previous
theorem, choosing $2^{\lambda}$ stationary sets
$W\subseteq S$ that are distinct modulo the cub ideal.
(In more detail, let $V$ and $W$ be two of these stationary sets.  Then
$\sm(M^W,\Mun^W)
\bigtriangleup
\sm(M^V,\Mun^V) \supseteq W \union V$.
Thus by Lemma~\ref{defrep},
$(M^W)\not\iso (M^V)$.)

\jtbnumpar{Remark}  If we add the requirement $\lambda^{\chi_1({\bf K})}=
\lambda$ we can deduce that $\lambda$ is {\bf K}-inaccessible from
Lemma~\ref{existinacc1}.
Applying Lemma~\ref{notsmoothstrat} we could omit DC2 from the
hypothesis list if $\kappa = \omega$.

\medskip
The assumption
in
Theorem \ref{maintheorem:smooth}
that {\bf K} is
$(<\lambda,<\lambda)$-bounded
is used only for the
construction of
the $M_{\alpha}$ for $\alpha \not\in S^{+}$.  We can weaken this model
theoretic hypothesis at the cost of strengthening the set theoretic
hypothesis.
We noted in Fact~\ref{justificationstrongboxa} that the set theoretic
hypotheses of the next theorem follow from V=L.

\begin{thm}
Suppose
$\Box^a_{\lambda,\kappa}(S)$ holds for an $S$
that does not reflect.
Then the hypothesis that {\bf K} is
$(<\lambda,<\lambda)$-bounded can be deleted from
Theorem~\ref{maintheorem:smooth}.
\end{thm}

\proof. The only use of this hypothesis is
the construction of $M_{\alpha}$ for $\alpha \not\in S^{+}$.
In this case we make our construction more uniform by demanding
for $\alpha \not \in S^{+}$ that
$M_{\alpha}$ is canonically prime over $\langle
M_{\beta}:\beta < \alpha\rangle$.
If $S$
does not reflect in $\alpha$ then there is a club $C\subseteq \alpha$
with $C
\inter S = \emptyset$.  By induction, for each $\delta \in C$, we have
$M_{\delta}$ canonically prime over $\langle M_{\beta}:\beta <
\delta\rangle$ .  Thus the chain $\langle M_{\delta}:\delta \in
C\rangle$ is {\bf K}-continuous
and we can choose $M_{\alpha}$ canonically
prime over it.
By Axiom Ch5,
$M_{\alpha}$ is canonically prime over $\langle
M_{\beta}:\beta < \alpha\rangle$ as required.

\medskip
Recall that
{\bf K} is $(<\lambda,[\mu,\lambda])$-bounded
if every chain with
cofinality between $\mu$ and $\lambda$ inclusive of
models that each have cardinality $< \lambda$ is bounded.
In a number of the examples we have adduced (\ref{variousexamples}),
{\bf K} is $(<\infty,[\mu,<\infty])$ bounded for appropriate $\mu$.
  Thus the model theoretic
hypothesis of the following theorem is reasonable.
The existence of stationary sets that do
not reflect in $\delta$ of small cofinality
is provable if {\bf V = L} and
is consistent with large cardinal hypotheses.


\begin{thm}
Fix $\kappa, \mu < \lambda$.
Suppose
$\Box^a_{\lambda,\kappa}(S)$ holds for some stationary subset $S$
of $\lambda$ that satisfies:
$$\hbox{\rm if $S$ reflects in $\delta$ then $\cf(\delta)\geq \mu.$}$$
Then the hypothesis that
{\bf K} is
$(<\lambda,<\lambda)$-bounded
can be replaced in
Theorem~\ref{maintheorem:smooth}
by assuming that
{\bf K} is $(<\lambda,[\mu,\lambda])$-bounded.
\end{thm}

\proof.  Again we must construct $M_{\alpha}$ for $\alpha \not \in S^+$.
If $S$ reflects in $\alpha$, $\cf(\alpha) \geq \mu$ so $\Mun\restrict
\alpha$ is bounded.  Since $\lambda$ is {\bf K}-inaccessible
(Definition~\ref{definacc})
we can
choose $M_{\alpha}$ to bound $\Mun\restrict \alpha$ and
with $|M_{\alpha}|< \lambda$.    If $S$ does not reflect in $\alpha$,
write $\alpha$ as a limit of ordinals $\beta_i$ of smaller cofinality.
By induction $M_{\beta_i}$
is canonically prime over $\Mun\restrict \beta_i$
and taking $M_{\alpha}$ canonically prime over the $M_{\beta_i}$
(using Axiom Ch5)
we
finish.

\subsection{Invariants modulo $\id(C^*)$}

In this subsection we show if {\bf K} is not smooth then for many
$\lambda$ we can code stationary subsets of $\lambda$ by assigning
invariants modulo the
ideal $\id(C^*)$ by the function $\som$.
We now replace the hypothesis that cpr behaves on
subseqences by assuming {\bf K} is weakly $(<\lambda,\theta)$-closed for
certain $\theta$; we use $\Box^b$ rather than $\Box^a$ but these have the
same set theoretic strength.

Again we first give
the technical  version of the main result with the parameters
and
reliance on the axioms enunciated in
Section~\ref{smaade} and Section~\ref{smatops}
stated explicitly.
We need to vary the meaning of the phrase, `a winning strategy against
$\Mun|C_{\alpha}$' by changing the game played on $C_{\alpha}$.

\jtbnot
Fix a square sequence $\langle C_{\alpha}:\alpha \in S\rangle$.
Suppose Player I has a winning strategy for Game 1$(\lambda,\kappa)$.
We modify our earlier notion of
what is meant by saying a certain $M_{\alpha}$ is chosen by playing
Player I's strategy on $\Mun|C_{\alpha}$ to a notion that is appropriate
for the proof of the next theorem.

Let $\langle c_{\beta}:\beta < \beta_0\rangle$ enumerate $C_{\alpha}$.
Denote $C_{\alpha} \cup \{\gamma +1:
\gamma \in \nacc[C_{\alpha}]\}$ by
$\hat C_{\alpha}$.
We attach to $\Mun|C_\alpha$ two sequences
$\langle L_i:i< \otp(C_{\alpha})\rangle$ and
$\langle P_i:i< \otp(C_{\alpha})\rangle$ by
setting
$$P_{\gamma} = M_{c_{\gamma}}$$
$$L_{\gamma} = M_{c_{\gamma} + 1} \text{ if $\gamma \in \nacc[C_{\alpha}]$ }$$
$$L_{\gamma} = M_{c_{\gamma}} \text{ if $\gamma \in \acc[C_{\alpha}]$ }$$
We say $M_{\beta}$ for $\beta = \alpha$ or $\beta \in C_{\alpha}$
is chosen by Player I's winning
strategy on $M|C_{\alpha}$ if the sequence
$\langle \Lun,\Pun\rangle$
associated with $C_{\alpha}\inter \beta \union \{\beta\}$ is
\begin{enumerate}
\item an initial segment of a play of Game 1 $(\lambda,\kappa)$ and
\item Player I's moves in this game follow his winning strategy.
\end{enumerate}

In defining this play of the game we have restrained Player II's
moves somewhat (as $P_{\gamma} = L_{\gamma}$ if $\gamma \in \acc[C_{\alpha}]$.
But this just makes it even easier for Player I to
play his winning strategy.
%

\begin{thm}
\label{maintheorem:smooth1}
Fix regular cardinals $\kappa, \theta,R, \lambda$ with
$\kappa \neq \theta$, $\theta < \lambda$,  $R\leq \lambda$,
and
$\chik, \kappa^+ < \lambda$.
Suppose the following.
\begin{enumerate}
\item
{\bf K} is an adequate class.
\item   Player I has a winning strategy for Game 1 $(\lambda,\kappa)$.
\item
$\lambda$ is a {\bf K}-inaccessible cardinal and for some
$S,S_1,S_2 \subseteq \lambda$
and {$C^{**}=\langle C_{\alpha}:\alpha \in S\rangle$,}
$\Box^b_{\lambda,\kappa,\theta,R}(S,S_1,S_2)$ is witnessed by $C^{**}$.
Let $S_0 = S_1 \union S_2$ and $C^*= C^{**}|S_0$.
\item
{\bf K} is $(<\lambda,[R,\lambda])$-bounded.
\item
 {\bf K} is $(<\lambda,\theta)$-closed.
\item
{\bf K} is $(<\lambda,\lambda)$-closed.
\end{enumerate}
Then
there is a model $M$ with $\som(M,C^*) = (S_2/\fil(C^*))$.
\end{thm}

\proof.
For each
$\alpha< \lambda$ we define
$M_{\alpha}$ as follows.
We will guarantee
$$\som(\Mun, M,C^*) = S_2 \text{ modulo } \fil(C^*).$$
The construction proceeds by induction.  There are a number of cases
depending on whether $\alpha\in S_1, S_2$ etc.
Choose $M_{\alpha}$ by the first of the following conditions that
applies to $\alpha$.
\begin{description}
\item[Case I.]
$\alpha \in (S_1
\cup \bigcup_{\delta\in S_1} \acc[C_{\delta}]
\cup \bigcup_{\delta\in S_1} \{\gamma+1:\gamma \in \nacc[C_{\delta}]
\}$:
Choose any $\beta\in S$ with $\alpha \in C_{\beta}$ (if $\alpha$ is in
$S$ let $\beta = \alpha$).
Apply
Player I's winning strategy for  Game 1 to $\Mun|C_{\beta}$ to
choose $M_{\alpha}$.  Again the coherence conditions in the definition
 of the square sequence guarantee that the particular
 choice of $\beta$ is
immaterial.  Note that for $\alpha$ not in $S_1$, the definition of
Player I winning the game guarantees that $M_{\alpha}$ is canonically
prime over $C_{\beta} \inter \alpha$.
\item[Case II.]
For any  successor ordinals not yet covered, say
$\beta = \gamma+1$, choose $M_{\beta}$
as a {\bf K}-extension of $M_{\alpha}$.
Thus in the rest of the cases we may assume $\alpha$ is a limit ordinal.
\item[Case III.]
$\alpha \in
\bigcup_{\delta\in S_0}\nacc[C_{\delta}]$:
Then $\cf(\alpha) = \theta$ so by the fifth hypothesis we may choose
$M_{\alpha} = \cup_{\beta<\alpha} M_{\beta}$ provided that $\Mun|\alpha$
is bounded.  If $\alpha \in S$, $\Mun|\alpha$ is bounded by the
canonically prime model over $\Mun|C_{\alpha}$ (which exists by the
argument for Case V).  If $\alpha \not \in S$ then
Definition~\ref{boxb}  guarantees $\cf(\alpha) \geq R$.
Choose $M_{\alpha}$
as a {\bf K}-extension of $M_{\beta}$ for each $\beta < \alpha$ by
$(<\lambda,\geq R)$-boundedness.
\item[Case IV.]
$\alpha \in S_2$:
By Case III
$\Mun | \nacc[C_{\alpha}]$
is {\bf K}-continuous; choose
$M_\alpha$ canonically prime over
$\Mun | \nacc[C_{\alpha}]$.
\item[Case V.]  All remaining ordinals
$\alpha \in S$:
Our construction guarantees that
$\Mun|C_{\alpha}$ is continuous as $C_{\alpha} \subseteq S - S_1$.
Choose
$M_{\alpha}$ canonically prime over $\Mun|C_{\alpha}$;
\item[Case VI.]
$\alpha$ is a limit ordinal and $\alpha \not \in S$:
Then Definition~\ref{boxb}  guarantees $\cf(\alpha) \geq R$.
Choose $M_{\alpha}$
as a {\bf K}-extension of $M_{\beta}$ for each $\beta < \alpha$ by
$(<\lambda,\geq R)$-boundedness.
\end{description}

Now we show that
$\som(\Mun,C^*,M)$
intersects the stationary set $S_0$
in $S_2$.  If $\alpha\in S_1$ the play of Game 1 guarantees that
there is no
$M'_{\alpha}$ such that
$\Mun|\nacc[C_\alpha]\union \{M'_{\alpha}\}$ is essentially {\bf K}
continuous.
Thus by condition iii) of Definition~\ref{defsom}
$\alpha\not \in
\som(\Mun,C^*,M)$
But if $\alpha \in S_2$, $\alpha$ is in neither $S_1$ nor
$\cup_{\beta\in S_0}\nacc[C_{\beta}]$ (since all elements of the
second set have cofinality $\theta$ and $\cf(\alpha) = \kappa$).
Thus, by Case IV  of the construction
$M_{\alpha}$ is canonically prime over $\Mun|\nacc[C_\alpha]$
and since {\bf K} is $(<\lambda,\lambda)$-closed $M_{\alpha} \subm M$.
Condition ii) in the definition of $\som$ is guaranteed since
$\Mun|C_{\alpha}$ is {\bf K}-continuous (as $C_\alpha \inter S_1 =
\emptyset$ and all points of non {\bf K}-continuity are in $S_1$).
Condition i), $M_{\gamma} = \cup_{\beta <\gamma}M_{\beta}$ for
$\gamma\in \nacc[C_{\alpha}]$, is guaranteed by Case III  of the
construction.

\jtbnumpar{Remark}
\label{weak}
At the cost of assuming that
{\bf K} is  $(<\lambda,<\kappa)$-smooth and $\theta<\kappa$,
hypothesis v) could be
weakened to
 {\bf K} is weakly $(<\lambda,\theta)$-closed.

Again we rephrase the result to emphasize the salient hypotheses.
\begin{thm}[ZFC]
\label{secondquotable}
Let {\bf K} be an
adequate class satisfying DC1 and DC2.
Fix $\lambda >\chik$ and $\theta$ such that
 \begin{enumerate}
\item
$\lambda$ is
{\bf K}-inaccessible,
\item $\lambda$  is
a succcessor of a regular cardinal and $\theta^+$ is less than $
\lambda$,
\item and {\bf K} is
 \begin{enumerate}
\item $(<\lambda,<\lambda)$-bounded,
\item
weakly $(<\lambda,\theta)$-closed,
\item $(<\lambda,\lambda)$-closed.
\end{enumerate}
\end{enumerate}
If {\bf K} has fewer than $2^{\lambda}$ models with cardinality
$\lambda$ and $\kappa^+<\lambda$
then {\bf K} is $(<\lambda,<\kappa)$-smooth.
\end{thm}

\proof.  Fix $\kappa$ with $\kappa^+<\lambda$
such that {\bf K}
is not $(<\lambda,\kappa)$-smooth.
By Lemma~\ref{mainnotsmoothstrat} Player I
has a winning strategy in Game 1 $(\lambda,\kappa)$.
By Theorem~\ref{justificationboxb}\,i),
$\Box_{\lambda,\kappa,\theta,\aleph_0}$ holds.  Now a very slight
variant of the proof of
\ref{maintheorem:smooth1} shows
there exist $2^{\lambda}$ models $M_i$
with the $\som(M_i)$ distinct modulo $\fil(C^*)$.
(Namely, since {\bf K} is $(<\lambda,<\lambda)$-bounded
we do not need to worry about the cofinality of $\alpha$ in
Case VI.)
\medskip

This shows that if {\bf K} is not smooth at some $\kappa$ then there
will be many models in power $\lambda$ for many $\lambda>\kappa$
satisfying certain model theoretic hypotheses.

If V=L we can waive the boundedness hypothesis.

\begin{thm}[V=L]
\label{nobound}
Let $\lambda > \chik$ and
not weakly compact be {\bf K}-inaccessible.
Suppose the adequate class {\bf K} satisfies DC1, DC2 and
is
\begin{enumerate}
\item
  weakly $(<\lambda,\mu)$-closed, for some $\mu \leq \lambda$,
\item $(<\lambda,\lambda)$-closed.
\end{enumerate}
If {\bf K} has fewer than $2^{\lambda}$ models in  power $\lambda$
and $\kappa^+ < \lambda$
then {\bf K} is $(<\lambda, \kappa)$ smooth.
\end{thm}

\proof.  Since V=L,
Lemma~\ref{justificationboxb} implies
$\Box^b_{\lambda,\kappa,\lambda,\mu}$ holds.
The result now follows from
Theorem~\ref{maintheorem:smooth1} taking
$\lambda$ as $\theta$ and $\mu$ as $R$.
We observed after Definition~\ref{defclosed} that
$(<\lambda,\lambda)$-closed
implies
$(<\lambda,\lambda)$-bounded.

\medskip
We have shown in this section that each of the variants of $\Box$
discussed in Section~\ref{smast} suffice to show that a nonsmooth {\bf
K} codes stationary subsets of $\lambda$ for many $\lambda$.

\section{The Monster Model}
\label{smamons}

We begin this section by recapitulating the assumptions that we will
make in developing the structure theory.  Then we show that under these
assumptions we can prove the existence of a monster model and prove the
equivalence between `homogenous-universal' and `saturated'.

\jtbnumpar{Convention}  $\chi = \chi_{\bf K}$ is a cardinal satisfying
Axioms S0 and S1 stated in Section~\ref{smaade}.
We assume $\mu \geq \chi$ throughout this section.

\jtbnot  $\kd$ denotes a compatibility class of {\bf K}.

\jtbdef We say {\bf K} satisfies the {\em joint embedding property} if
for any $A, B\in {\bf K}$ there is a $C\in {\bf K}$ and
{\bf K}-embeddings of $A$ and $B$ into $C$.

At this point we extend the axioms from Section~\ref{smaade} by adding the
smoothness hypothesis we have justified in the last few sections.
We have shown that under reasonable set theoretic hypotheses the failure
of these smoothness conditions allows us to code stationary subsets of
$\lambda$ for a proper class of $\lambda$.

\jtbnumpar{Assumptions}
\label{newadequate}
{\bf
We assume in this section axiom groups A and C,
Axioms D1 and D2 from group D,
(all from \cite{BaldwinShelahprimali}),
axioms
Ch1', Ch2', Ch4,
L1,
S0 and S1 from
Section~\ref{smaade} and the following smoothness conditions.
{\bf K} is $(<\infty,\geq \chin)$ smooth and $(<\infty,> \chin)$
fully smooth.
Thus, we assume {\bf K} is $(<\infty,\infty)$-smooth.
Finally we assume that {\bf K} satisfies the joint embedding property.
We call such a class {\em fully adequate}.
}

The assumption of the joint embedding property is purely a notational
convenience.  We have just restricted from {\bf K} to a single
compatibility class in
${\bf K}$.  Thus, the notions that in
\cite{Shelahuniversal}
are written, e.g., $({\Dscr},\mu)$-homogeneous here
become $({\bf K},\mu)$-homogeneous with no loss in generality.  We could
in fact drop the {\bf K} altogether.

This class is called fully adequate because (modulo $V=L$, see
Conclusion~\ref{concl}\,ii)) any fully adequate class either has a
unique homogeneous-universal model or many models.
We did not include the axiom  Ch5, `cpr behaves on subsequences',
since we can rely on  Theorem~\ref{maintheorem:smooth1} to obtain
smoothness
without that hypothesis.  We did assume that {\bf K} is $(<\lambda,\mu)$
closed for large $\mu$ so the other hypotheses of that theorem are
fulfilled.

The following observation allows us to perform the required
constructions.

\begin{lemm}
\label{smallprime}
If {\bf K} satisfies the $\lambda$-L\"{o}wenheim-Skolem
property, $\lambda \geq \chin$,
$\langle M_i: i < \kappa\rangle$
is a
chain of models inside $M$ with each $|M_i| < \lambda$ and $\cf(\kappa)<
\lambda$ then there is a  canonically prime model $M'$ over
$\langle M_i: i < \kappa\rangle$
with $|M'| < \lambda$.
\end{lemm}
\proof.  If $\cf(\kappa) \geq \chik$, $M' = \cup_{i <\kappa} M_i$ is
the required model.  If not, note that by the
$\lambda$-L\"{o}wenheim-Skolem
property there is an $N \subm M$ containing the union with $|N| <
\lambda$. By $(<\infty,\geq \chik)$ smoothness any canonically prime
over
$\langle M_i: i < \kappa\rangle$
can be embedded in $N$.


\jtbdef
The model $M$ is $({\bf K},\mu)$-homogeneous if
\begin{enumerate}
\item
for every $N_0
\subm M$ and every $N_1$ with $N_0\subm N_1 \in {\bf K}$ and $|N_1| \leq
\mu$ there is a {\bf K}-embedding of $N_1$ into $M$ over $N_0$;
\item every $N\in {\bf K}$ with cardinality less than $\mu$ can be {\bf
K}-embedded in $M$.
\end{enumerate}

There is a certain entymological sense in labeling this notion a kind of
saturation.  The argument for homogeneity is that naming {\bf K}
describes the level of universality and we need only indicate the
homogeneity again.
In any case Shelah established this convention
almost twenty years ago in \cite{Shelahfindiag}.
We identify this algebraic notion with a
type realization notion in Theorem~\ref{sathom}.

\begin{lemm}
Suppose {\bf K} is $(<\mu,<\mu)$-bounded and $(<\mu,\mu)$-closed.
\begin{enumerate}
\item
If $\mu
\geq
\chik$,
is regular,
{\bf K}-inaccessible and satisfies $\mu^{<\mu} = \mu$,
there is a
$({\bf K},\mu)$-homogeneous model of power $\mu$.
\item
If, in addition {\bf K} is $(<\mu,\leq \mu)$ smooth this
$({\bf K},\mu)$-homogeneous model is unique up to isomorphism.
\end{enumerate}
\end{lemm}

\Proof.  We define an increasing chain $\langle M_i: i < \mu \rangle$ by
induction; the union of the $M_i$ is the required model.  Let $M_0$ be
any element of ${\bf K}_{<\mu}$.
Fix an enumeration $\langle N_{\beta}:\beta < \mu\rangle$ of all
isomorphism types in ${\bf K}_{<\mu}$.  There are
only $\mu$ such since $\mu^{<\mu} =\mu$.
Given $M_i$, with  $|M_i| < \mu$ we
define $M_{i+1}$ as a bound for a sequence $M_{i,j}$ with
$j <
|i| + |M_i|^{|M_i|} =
\alpha < \mu$.  First let
$G_j = \langle A_j, B_j, f_j\rangle$ for $j <
\alpha$
be a list of
all triples such that $f_j$ is an isomorphism of $A_j$ onto a
{\bf K}-submodel of
$M_i$ and $A_j \subm B_j$ and $B_j \iso N_{\beta}$ for some
$\beta < i$.
(Note that the $B_j$ are specified only up to isomorphism; a given
isomorphism type of $A_j$ will occur many times in the list depending on
various embeddings $f_j$ into $M_i$.)
Now, $M_{i,0} = M_{i}$;  $M_{i,j+1}$ is the amalgam of $M_{i,j}$ and
$B_j$ over $A_j$ (via $f_j$ and the identity map).  If $\delta$ is a
limit ordinal less than $\alpha$, $M_{i,\delta}$ is any bound of
$\langle M_{i,j}:j <\delta \rangle$ with $|M_{i,\delta}| < \mu$.
$M_{i+1}$ is a bound for the $M_{i,j}$.  By regularity of $\mu$
for limit $\delta < \mu$,
each
$M_{i,\delta}$ and $M_{\delta}$ have cardinality less than $\mu$.

It is easy to see that $M$ is homogeneous since if $f:N_0 \mapsto M $ is
a {\bf K}-embedding and $N_0 \subm N_1$ with $|N_1| < \mu$, f was
extended to a map into some $M_{i,j}$ at some stage in the construction
and $M_{i,j}\subm M$.

The uniqueness of the
$({\bf K}, \mu)$-homogeneous model
now follows by
the usual back and forth argument to show any two
$({\bf K}, \mu)$-homogeneous models $M$ and $N$
of power $\mu$ are isomorphic.
But smoothness is crucial.  At a limit
stage $\delta$,
one takes the canonically prime model $M_{\delta}$
over an initial segment of the
sequence of submodels of $M$ and embeds it as a submodel $N_{\delta}$
of $N$.  In order to continue the induction we must
know
$M_{\delta}$
is a strong submodel of $M$ and this
is guaranteed by smoothness.

\jtbnumpar{Conclusion}
\label{concl}
\begin{enumerate}
\item
For any fully adequate
${\bf K}$ that is $(<\infty,<\infty)$-bounded
there is (in some cardinal $\mu$) a
unique
$({\bf K},\mu)$-homogeneous model.
\item If $V=L$, we can omit the boundedness hypothesis
(by Theorem~\ref{nobound}.
\item
We will call
unique
$({\bf K},\mu)$-homogeneous model,
$\Mscr$, the monster model.  From now on all sets and
models are contained in $\Mscr$.
\end{enumerate}

\jtbnumpar{Remark}
  This formalism encompasses the constructions by Hrushovski
\cite{Hrustableplane}
of $\aleph_0$-categorical stable pseudoplanes.  An
underlying (but unexpressed)
theme of his
constructions is to generalize the Fraiss\'e-J\'onsson construction by
a weakening of homogeneity.  He does not demand that any isomorphism of
finite substructures extend to an automorphism but only an isomorphism
of submodels that are `strong substructures' (where strong varies
slightly with the construction).  This is exactly encapsuled in the
formalism here.
This viewpoint is pursued in \cite{Baldwinnewfam}.
Of course in Hrushovski's case the real point is the
delicate proof of amalgamation and $\omega$ is trivially {\bf
K}-inaccessible.  We assume amalgamation and worry about  inacessibility
and smoothness in larger cardinals.
\jtbdef
\begin{enumerate}
\item
The {\em type} of $\abar$ over $A$ (for $\abar, A \subseteq
\Mscr$) is the orbit of $\abar$ under the automorphisms of $\Mscr$ that
fix $A$ pointwise.
We write $p = \tp(\abar;A)$ for this orbit.
\item
$p$ is a $k$-type if $\lg(\abar) = k$.
\item
The {\em type} of $B$ over $A$ (for $B, A \subseteq \Mscr$) is the
type of some (fixed) enumeration of $B$.
\item $S^{k}(A)$ denotes the collection of all $k$-types over $A$.
\end{enumerate}
We will often write $p$, $q$, etc. for types.  This notion is really of
interest only when $\lg (\abar) \leq \mu$; despite the suggestive
notation, $k$-type, we may deal with types of infinite length.
We will write $S(A)$ to mean $S^k(A)$ for some $k <\mu$ whose exact
identity is not important at the moment.

\jtbdef
\begin{enumerate}
\item The type $p \in S(A)$ is  {\em realized} by $\cbar\in N$ with $A
\subseteq N \subm \Mscr$ if $\cbar$ is member of the orbit $p$.
\item $N \subm \Mscr$ is
$({\bf K}, \lambda)$-{\em saturated} if for every
$M\subm N$ with $|M| < \lambda$, every 1-type over $M$ is realized in
$N$.
\end{enumerate}

\begin{thm}
\label{sathom}
Let $\lambda \geq \chik$ be {\bf K}-inaccessible.  Then $M$
is $({\bf K}, \lambda)$-saturated if and only if $M$ is
$({\bf K},\lambda)$-homogeneous.
\end{thm}

The proof follows that of Proposition 2.4 of \cite{ShelahMakkaicat} line
for line with one exception.  If we consider those stages $\delta$ in
the construction
where
$\cf(\delta) < \chin$, we cannot
form $M_{\delta}$
just by taking unions.  However, any
canonically prime model over the initial segment of the construction
will work by smoothness and Lemma~\ref{smallprime}.

\section{Problems}
\label{conc}
\begin{quest}  Can one give more precise information on the class of
cardinals in which an adequate class {\bf K} has a model.
\end{quest}

In \ref{reinonsmooth} we gave a definition of $\subm$ on the class of
$\aleph_1$-saturated models of the theory $T =\REI_{\omega}$ under which
this class is not $(<\aleph_1,\omega)$-smooth.  Thus, by our main result
$T$ has the
maximum number of $\aleph_1$-saturated models in power $\lambda$ (if
e.g. $\lambda = \mu^+$).

\begin{quest}  Define $\leq$ on the class of $\aleph_1$-saturated models
of an strictly stable  with {\em didop} \cite{Shelahmaingap1}
(or perhaps if {\bf K} is not finitely controlled in the sense of
\cite{Hrufinbase})
so the
class is not smooth.
\end{quest}

There are strictly stable theories with fewer than the
maximal number of $\aleph_1$-saturated models in most $\lambda$.
See: Example 8 page 8 of \cite{Baldwinbook}.

\begin{quest}  Formalize the notion of coding a stationary set to
encompass the examples we have described and clarify the distinctions
described at the beginning of Section~\ref{smainv}.
\end{quest}

We have developed this paper entirely in the context of cpr models.  In
a forthcoming work we replace this fundamental concept by axioms for
winning games
similar to Game 1 $(\lambda,\kappa)$ and establish smoothness in that
context.  The cost is stronger set theory (but V=L suffices).
\bibliography{ssgroups}
\bibliographystyle{plain}
\end{document}